\documentclass[journal,twoside,web]{ieeecolor}

\usepackage{generic}
\usepackage{cite}
\usepackage{amsmath,amssymb,amsfonts}
\usepackage{graphicx}
\usepackage{textcomp}
\usepackage{epstopdf}
\usepackage{graphics} % for pdf, bitmapped graphics files
\usepackage{epsfig} % for postscript graphics files
\usepackage{amsmath} % assumes amsmath package installed
\usepackage{amssymb,cite} % assumes amsmath package installed
\usepackage[dvipsnames]{xcolor}
\usepackage{bbm}
\usepackage{enumerate}
\usepackage{algorithm}
\usepackage{amstext}
\usepackage{subfig}
\usepackage{setspace}
\usepackage{todonotes}
\allowdisplaybreaks % to allow breaking align blocks and others over two columns/pages
\usepackage[noend]{algpseudocode}
\newcommand{\RN}[1]{%
	\textup{\uppercase\expandafter{\romannumeral#1}}%
}
\makeatletter
\def\BState{\State\hskip-\ALG@thistlm}
\makeatother
\graphicspath{{figures/}}
{
	\newtheorem{assumption}{Assumption}
}

{
	\newtheorem{lemma}{Lemma}
}

{
	\newtheorem{theorem}{Theorem}
}
{
	\newtheorem{remark}{Remark}
}
{
	\newtheorem{proposition}{Proposition}
}
{
	
}
{
	\newtheorem{corollary}{Corollary}
}
\def\BibTeX{{\rm B\kern-.05em{\sc i\kern-.025em b}\kern-.08em
  T\kern-.1667em\lower.7ex\hbox{E}\kern-.125emX}}
\begin{document}
\title{Asynchronous Distributed Optimization
	via Dual Decomposition and Block Coordinate Subgradient Methods}
\author{Yankai Lin, \IEEEmembership{Student Member, IEEE}, Iman Shames, \IEEEmembership{Member, IEEE} and Dragan Ne\v{s}i\'{c}, \IEEEmembership{Fellow, IEEE}
\thanks{This work was supported by the Australian Research Council under the Discovery Project DP170104099.}
\thanks{The first and the third authors are with the Department of Electrical and Electronic Engineering, The University of Melbourne, Parkville, 3010, Victoria, Australia (email: yankail@student.unimelb.edu.au; dnesic@unimelb.edu.au). The second author is with the School of Engineering, The Australian National University, Acton, 0200, ACT, Austrlia (email:iman.shames@anu.edu.au)}}
%\thanks{S. B. Author, Jr., was with Rice University, Houston, TX 77005 USA. He is 
%now with the Department of Physics, Colorado State University, Fort Collins, 
%CO 80523 USA (e-mail: author@lamar.colostate.edu).}
%\thanks{T. C. Author is with 
%the Electrical Engineering Department, University of Colorado, Boulder, CO 
%80309 USA, on leave from the National Research Institute for Metals, 
%Tsukuba, Japan (e-mail: author@nrim.go.jp).}}

\maketitle

\begin{abstract}
We study the problem of minimizing the sum of
potentially non-differentiable convex cost functions with partially overlapping dependences in an asynchronous manner, where communication in the network is not coordinated. We study the behavior of an asynchronous algorithm based on dual decomposition and block coordinate subgradient methods under assumptions weaker than those used in the literature. At the same time, we allow different agents to use local stepsizes with no global coordination. Sufficient conditions are provided for almost sure convergence to the solution of the optimization problem. Under additional assumptions, we establish a sublinear convergence rate that in turn can be strengthened to linear convergence rate if the problem is strongly convex and has Lipschitz gradients. We also extend available results in the literature by allowing multiple and potentially overlapping blocks to be updated at the same time with non-uniform and potentially time varying probabilities assigned to different blocks. A numerical example is provided to illustrate the effectiveness of the algorithm.
\end{abstract}


\section{INTRODUCTION}

% \IEEEPARstart{D}{istributed} optimization has received significant attention due to its wide applications in areas such as resource allocation \cite{XJB} (e.g.~power networks, communication networks and so on), distributed model predictive control (MPC) \cite{MN}, and many others. These problems typically involve a large number of agents that are interconnected via a network without a centralized coordinator.

%Specifically, we allow communications between agents to be asynchronous and the asynchrony is modelled by assigning random timers to agents in the network dictating whether an update is happening or not. These random timers account for availability of access to the network that may be the result of capacity and bandwidth constraints of the network or potential link failures.

\IEEEPARstart{D}{istributed} optimization algorithms can be roughly classified as primal and dual methods. For primal methods, there is a large body of literature built upon the seminal work by Bertsekas and Tsitsiklis \cite{BT} where numerical computations are assigned to a group of processors (or agents) that exchange information with their neighbors. 
%A central tool to the design of distributed algorithms for minimizing the sum of the costs is consensus theory as it allows agents to solve their own local problems independently while guaranteeing that all agents converge to the same value. 
In \cite{NO}, a distributed (sub)gradient method was proposed, which can be viewed as a first order primal method based on consensus theory. 
%At each step, every agent updates its estimate of the solution via a linear combination of a gradient descent step and a consensus step. Under the assumption of bounded gradient and diminishing stepsizes, it is proved that distributed gradient descent converges exactly to the solution of the problem. 
The consensus based algorithm is later extended to other settings including algorithms using constant stepsizes \cite{SLWY}, stochastic networks \cite{XZSX}, and algorithms with non-uniform stepsizes \cite{LRYG}. The existing dual methods are based on the dual decomposition idea. In \cite{TTM}, a synchronous distributed algorithm for convex problems over a fixed network with bounded communication delays is proposed. The case of time-varying networks is considered in \cite{FMGP}. This class of methods also include the distributed augmented Lagrangian method \cite{FGGN} and the alternating direction method of multipliers (ADMM) \cite{BPCPE,TGSJ,WO}.

\paragraph{The Main Contributions}
First, we study an algorithm based on the dual decomposition where agents solve a primal problem locally and independently. Then, they communicate these solutions to a subset of their neighbors whenever they finish the necessary computations and if communication is possible. Particularly, we consider solving optimization problems in an asynchronous environment where the notion of asynchrony is adopted from \cite[Chapter 1.4]{BT}. Asynchrony refers to both independent computation and asynchronous communication. We investigate the conditions under which the algorithm converges under such asynchronous conditions. In this work, we propose an asynchronous dual supergradient algorithm and give sufficient conditions on the network and stepsizes to guarantee almost sure convergence of the algorithm. As our algorithm is based on dual decomposition, differentiability of the primal problem does not play a role in our analysis. 

Second, under mild assumptions, we allow the use of local stepsizes where each agent chooses and updates its own stepsizes. Moreover, we assume that each component is updated infinitely often. This includes commonly used deterministic updating rules such as round robin and persistently exciting \cite{TN} in the same unifying framework. Since we assume each block is updated infinitely often, our main result still holds if multiple and overlapping blocks are updating at the same time. 

Lastly, we prove our main convergence result using tools from asynchronous stochastic approximation theory \cite{Borkar1,RBQ} which does not rely on strong convexity and Lipschitzness of the gradient of the cost function. To this end, we extend the results in \cite[Chapter 7]{Borkar1} to the non-smooth case which we believe is an important result in itself. Consequently, we prove our main result for the general case without assuming strict or strong convexity of the problem. This, inevitably, results in our convergence rate estimates to be slower than the results that consider the problem under more stringent conditions (e.g., the aforementioned strong convexity and Lipschitz conditions). However, it can be shown that, if stronger assumptions are imposed, we recover the linear convergence rate of the type shown in \cite{Wright} as a special case. 

%First, we assume only convexity of the cost function without appealing to its differentiability. This is in contrast with most existing literature that relies on assumptions such as strongly convex cost function with Lipschitz continuous or bounded gradients \cite{NO,NCN,XZSX}. Second, compared to the block subgradient algorithm in \cite{LSN} and \cite{CFNN} where it is assumed that each component has a positive probability of being updated at each time step, we only assume that each component is updated infinitely often. This allows us to include some commonly used deterministic updating rules such as round robin and persistently exciting \cite{TN} which are not included in the aforementioned references. Additionally, overlapping blocks are allowed in our analysis which is different from the set-up considered in \cite{NCN}, where each component is activated by a random local clock that is assumed to be independent and identically distributed (i.i.d.). In contrast to \cite{LSN} and \cite{CFNN}, where the stepsizes are assumed to be global in the sense that all edges that are updating use the same stepsize at each time step, we give sufficient conditions to allow the use of local stepsizes that may be different at each time step for different components. Lastly, our analysis is based on tools from asynchronous stochastic approximation theory \cite{Borkar1,RBQ}, and we allow the use of different local stepsize sequences for different coomponents compared to \cite{Borkar1,RBQ}. 
\paragraph{Comparison to the Literature}
For a recent survey paper on distributed optimization algorithms we refer the reader to \cite{Survey}. The most relevant ones are those based on dual methods and asynchronous algorithms. In the consensus based primal gradient algorithm considered in \cite{XZSX}, the asynchrony is captured by the random communication graph that is drawn i.i.d.~from a set of graphs with proper weights to ensure almost sure convergence. Heterogeneous constant stepsizes are used in the asynchronous network under the assumption that the cost functions have Lipschitz gradients. If the cost functions are strongly convex, linear convergence rate can be established. 

Other than gradient methods, augmented Lagrangian methods are also used in asynchronous distributed optimization \cite{FGGN}. It is assumed that each agent updates at least once in a finite number of steps \cite[Assumption 2]{FGGN}. Another relevant work is \cite{BHI} on asynchronous primal dual algorithms. The authors in \cite{BHI} developed a stochastic primal-dual algorithm based on stochastic coordinate descent. Almost sure convergence of the proposed algorithm is established under the assumption of i.i.d.~updates. In \cite{TTM}, the authors propose a dual algorithm for solving the optimal consensus algorithm in a synchronous fashion with diminishing stepsizes. In \cite{NCN}, the authors consider the dual algorithm of \cite{TTM} in an asynchronous network for strongly convex functions with Lipschitz gradients. A dual block coordinate ascent algorithm is proposed with linear convergence rate under the assumption of i.i.d.~updates.

% The major difference between our work and \cite{XZSX} is that, we take the approach of dual decomposition and the algorithm is proposed regarding the dual variables. Moreover, instead of updating the agents in an i.i.d.~fashion as done in \cite{XZSX}, we allow agents to have their own local clocks and update infinitely often which covers some important deterministic cases and cases where the probability distribution of selection may be time-varying (due to adversarial reasons for instance). Our assumption of infinitely often update covers both the stochastic case considered in \cite{XZSX} and the deterministic case considered in \cite{FGGN}. Moreover, \cite[Assumption 3]{FGGN} requires that only one agent can update at each time instant, which is not imposed in our work.
% In contrast, we extend the dual algorithm in \cite{TTM} to the asynchronous setting and allow agents to use local and uncoordinated stepsizes to achieve almost sure convergence of the dual variable. 

In this work, we relax the i.i.d.~assumption of the aforementioned works to infinitely often updates and allow agents to have their own local clocks and uncoordinated stepsizes. Consequently, we provide a framework that includes the deterministic update rules of the type presented in \cite{FGGN} and the stochastic update rules used in \cite{XZSX,LSN,CFNN} as special cases. Moreover, \cite[Assumption 3]{FGGN} requires that only one agent can update at each time instant, which is relaxed in our work. Additionally, unlike \cite{TTM} and \cite{NCN} where strict and strong convexity are assumed to ensure the dual function is differentiable and has a Lipschitz continuous gradient respectively, we prove almost sure convergence without appealing to the differentiability of the dual function. As a result, our dual convergence result holds for general convex primal cost functions, and we recover the result in \cite{NCN} if strong convexity and i.i.d.~updates are assumed. Lastly, compared to \cite{LSN}, we allow the use of uncoordinated local stepsizes where each agent updates its own stepsizes under mild assumptions. Moreover, we only assume that each component is updated infinitely often while in \cite{LSN} it is assumed that each block component of the subgradient vector is computed with a positive probability at each step.

% \paragraph{Organization}In Section~\ref{S-Notation}, we cover notation and some basic definitions. The problem is formulated in Section~\ref{S-PF}. In Section~\ref{sec3}, we introduce the asynchronous block coordinate ascent algorithm for the dual optimization problem and prove its convergence. Additional results under stronger conditions are given in Section~\ref{S-AR}. Numerical examples are presented in~\ref{S-NE} and are followed by conclusions and possible future research in Section~\ref{S-FW}. 
\section{NOTATION AND BASIC DEFINITIONS}\label{S-Notation}
Let $\mathbb{R}$ be the set of real numbers and $\mathbb{R}^{n}$ be the $n$-dimensional Euclidean space, $\mathbb{R}_{\geq 0}$ (resp. $\mathbb{R}_{>0}$) be the set of non-negative (resp. positive) real numbers. For $x,y\in\mathbb{R}^n$, $\langle x,y\rangle$ denotes the inner product in $\mathbb{R}^n$. The set of non-negative (resp. positive) integers are denoted by $\mathbb{Z}_{\geq 0}$ (resp. $\mathbb{Z}_{>0}$) and the finite set of non-negative integers from $0$ to $N\in\mathbb{Z}_{\geq 0}$ is denoted by $\mathbb{Z}_{N}$. Furthermore, $|\cdot|$ denotes the Euclidean norm of a vector $x\in \mathbb{R}^{n}$ and the corresponding induced norm of a matrix, respectively. The matrix $I_n$ is used to denote the $n$-dimensional identity matrix and $n$ will be dropped when the dimension is clear. The vector in $\mathbb{R}^N$ consisting of all ones is denoted by $\textbf{1}_N$. Let $C$ be a set, $|C|$ denotes its cardinality. The set $\mathbb{B}$ is the closed unit ball in $\mathbb{R}^n$ centred at the origin. For a closed set $S\subset\mathbb{R}^n$ and $\varepsilon>0$, $S+\varepsilon\mathbb{B}$ is the set $\{x\in\mathbb{R}^n:|x|_S\leq\varepsilon\}$, where $|x|_S:=\inf_{y\in S}|x-y|$. Set-valued mappings are indicated by the symbol $\rightrightarrows$. For a convex function $f(x):\mathbb{R}^n\rightarrow \mathbb{R}$, a \textit{subgradient} of $f(x)$ at $\bar{x}$ is a vector denoted by $g(\bar{x})\in\mathbb{R}^n$ such that $f(x)\geq f(\bar{x})+\langle g(\bar{x}),(x-\bar{x})\rangle$ for all $x\in\mathbb{R}^n$. The \textit{subdifferential} of $f(x)$ at $\bar{x}$ denoted by $\partial f(\bar{x}):\mathbb{R}^n\rightrightarrows\mathbb{R}^n$ is the set of all subgradients of $f(x)$ at $\bar{x}$. For a concave function $h(x):\mathbb{R}^n\rightarrow \mathbb{R}$, a vector $g(\bar{x})\in\mathbb{R}^n$ is a \textit{supergradient} of $h(x)$ at $\bar{x}$ if $-g(\bar{x})$ is a subgradient of $-h(x)$. The \textit{superdifferential} of $h(x)$ at $\bar{x}$ denoted by $\partial^+ h(\bar{x}):\mathbb{R}^n\rightrightarrows\mathbb{R}^n$ is the set of all supergradients of $h(x)$ at $\bar{x}$. All random variables considered in this paper are defined on a probability space $(\Omega,\mathcal{F},\mathbb{P})$, where $\Omega$ is the sample space, the $\sigma$-field $\mathcal{F}$ is the set of events and $\mathbb{P}$ is the probability measure defined on $(\Omega,\mathcal{F})$. For a vector sequence $x\in \mathbb{R}^n$, $x^+$ denotes its value at the next time instant, $x_{(i)}$ denotes the $i^{th}$ component of $x$, $x_{[i]}$ denotes the value of $x$ at time $i$, and $x_{(i),[j]}$ denotes the value of the $i^{th}$ component of $x$ at time $j$.
\section{PROBLEM FORMULATION}\label{S-PF}
We consider an optimization problem, where $N$ agents are used to cooperatively minimize the sum of $N$ individual cost functions that have partially overlapping dependences. The interconnection of the agents is described by an undirected graph $\mathcal{G}=(\mathcal{V},\mathcal{E})$, where $\mathcal{V}=\{1,2,...,N\}$ is the set of nodes, $\mathcal{E}\subseteq \mathcal{V}\times\mathcal{V}$ is the set of edges. An edge $(i,j)\in\mathcal{E}$ means that agent $i$ is able to access information from agent $j$. The (primal) problem we study in the work is given below
\begin{align}
 \underset{x\in \mathcal{X}\subseteq\mathbb{R}^{n}}{\text{min}}\quad & F(x):=\sum_{i=1}^{N}F_i(x_{(i)}) \label{PP} \\
\text{s.t. } \quad & E_{ij}x_{(i)}-E_{ji}x_{(j)}=0,\ \forall (i,j)\in\mathcal{E}, \notag
\end{align}
where $x_{(i)}\in\mathbb{R}^{n_i}$ is the local variable held by agent $i\in\mathcal{V}$, $x=(x_{(1)},x_{(2)},\ldots,x_{(N)}):=[x_{(1)}^T\ x_{(2)}^T\ ...\ x_{(N)}^T]^T\in\mathbb{R}^{n}$, $\sum_{i=1}^{N}n_i=n$ is the collection of all local variables, $F_i:\mathbb{R}^{n_i}\rightarrow\mathbb{R}$ is the individual cost function of agent $i$ and the matrices $E_{ij}$ and $E_{ji}$ are selection matrices for agents $i,j$ that characterize how the individual functions share common variables. The constraint set $\mathcal{X}=\mathcal{X}_1\times\cdots\times\mathcal{X}_N$ is assumed to be convex. Compared to \cite{TTM}, we allow $\mathcal{X}$ to be a non-strict subset of $\mathbb{R}^n$ and non-compact. Let $x^*\in \mathcal{X}$ denote a solution to problem (\ref{PP}) and $F^*:=F(x^*)$ be the optimal cost.
\begin{remark}\label{motiex}
	The problem (\ref{PP}) provides a framework for addressing a range of different problems and it can be shown that it is equivalent to the problem considered in \cite{AS}. Consider the optimal consensus problem of minimizing the sum of convex cost functions $\sum_{i=1}^{N}F_i(x)$, where $x\in\mathbb{R}^n$ is the decision variable. It can be equivalently re-formulated as
	\begin{align}\label{OCP}
	 \underset{\bar{x}\in \mathbb{R}^{nN}}{\text{min}}\quad & F(\bar{x}):=\sum_{i=1}^{N}F_i(x_{(i)})\\
	\text{s.t. } \quad & x_{(i)}-x_{(j)}=0,\ \forall i,j\in\mathcal{V}. \notag
	\end{align}
	It is clear that (\ref{OCP}) is equivalent to (\ref{PP}) with all of the selection matrices $E_{ij}$ being identity and $\mathcal{X}=\mathbb{R}^{nN}$. Moreover, some other problems such as distributed MPC, network flow optimization and resource allocation can also be re-formulated in the form of (\ref{PP}), e.g.~see \cite{NCN,AS}. 
	\hfill $\blacksquare$
\end{remark}
\begin{remark}\label{revise1}
% 	To prevent unnecessary communications, agents that do not share any common variables are not required to communicate. Therefore, 
% 	There are no communication links between agents $i$ and $j$ such that $(i,j)\notin\mathcal{E}$. Consequently, the communication graph $\mathcal{G}$ may not be complete or connected. The same partially overlapping set-up is also adopted in \cite{HNS}.
  Our problem formulation is similar to the set-up in \cite{HNS}.
% 	However, if the problem is fully coupled such as (\ref{OCP}), then the graph should be connected to ensure convergence. I
However, if the goal is to optimize a global decision variable then the graph should be connected to ensure convergence. 
	\hfill $\blacksquare$
\end{remark}

We now introduce the dual variables $\lambda_{(ij)}$, for agents $i,j\in\mathcal{V}$ that share common variables, which is associated with the constraint $E_{ij}x_{(i)}=E_{ji}x_{(j)}$. As a result, we use an edge-based communication structure by assigning the task of updating $\lambda_{(ij)}$ corresponding to the constraint $E_{ij}x_{(i)}-E_{ji}x_{(j)}=0$ and edge $(i,j)$ to agent $i$ if $(i,j)\in\overline{\mathcal{E}}:=\{(i,j)\in\mathcal{E}:i<j\}$ and to agent $j$ otherwise. Additionally, we denote the set of neighbors of agent $i$ that are not responsible for updating the corresponding dual vectors by $\mathcal{N}_i:=\{j\in\mathcal{V}:(i,j)\in\overline{\mathcal{E}}\}$. % and the set of the rest of neighbors of agent $i$ is denoted by $\mathcal{N}_i^-:=\{j\in\mathcal{V}:(j,i)\in\overline{\mathcal{E}}\}$. 
 Moreover, $\mathcal{I}:=\{i\in\mathcal{V}:\mathcal{N}_i\neq\emptyset\}$ denotes the set of agents that are responsible for updating at least one dual vector \footnote{Depending on different ways of writing the equality constraint in (\ref{PP}), there are multiple ways of assigning the task of updating dual variable which leads to different communication graphs and different constructions of $\overline{\mathcal{E}}$. However, these different dual vector assignments result in problems mathematically equivalent to the one discussed in our paper and can be analysed using exactly the same approach used in our paper {\it mutatis mutandis}.}.
% We define $\overline{\mathcal{N}}_i^+$ as an ordered set of all elements of $\mathcal{N}_i$ where the elements are ordered in an ascending order.		
  Hence, agent $i\in\mathcal{I}$ has access to $(\ldots,\lambda_{(ij)},\ldots)\in\mathbb{R}^{\bar{n}_i},j\in\mathcal{N}_i$ and the full dual vector is given by: $\lambda=(\ldots,\lambda_{(ij)},\ldots)\in\mathbb{R}^{\bar{n}}$ with $(i,j)\in\overline{\mathcal{E}}$. By definition, $\bar{n}_i=\sum_{j}\bar{n}_{ij}:j\in\mathcal{N}_i$ and $\bar{n}=\sum_{i}\bar{n}_{i}:i\in\mathcal{I}$. Here, $\bar{n}_{ij}$ is the dimension of the dual variable $\lambda_{(ij)}$ corresponding to the edge $(i,j)\in\overline{\mathcal{E}}$, $\bar{n}_i$ is the overall dimension of the dual variables assigned to agent $i\in\mathcal{I}$, and $\bar{n}$ is the dimension of the dual vector $\lambda$. It can be seen that the number of dual vectors is equal to $|\overline{\mathcal{E}}|$. The Lagrangian $L$ is given by
\begin{equation}\label{lagrangian}
L(x, \lambda):=\sum_{i=1}^{N}F_i(x_{(i)})+\sum_{(i,j)\in\overline{\mathcal{E}}}\langle\lambda_{(ij)},E_{ij}x_{(i)}-E_{ji}x_{(j)}\rangle.
\end{equation}
The corresponding dual function is given by
\begin{equation}\label{dual}
Q(\lambda):=\text{inf }L(x, \lambda), \quad \text{s.t. }\quad x\in\mathcal{X}.
\end{equation}
Due to the finite sum structure of $F(x)$, (\ref{dual}) can be rewritten as $Q(\lambda)=\sum_{i=1}^{N}q_i(\lambda)$, where
\begin{equation}\label{dualind}
q_i(\lambda)=\text{inf } l_i(x_{(i)},\lambda), \quad \text{s.t. } \quad x_{(i)}\in \mathcal{X}_i,
\end{equation}
with $l_i(x_{(i)},\lambda):=F_i(x_{(i)})+\sum_{j:(i,j)\in\overline{\mathcal{E}}}\langle\lambda_{(ij)},E_{ij}x_{(i)}\rangle-\sum_{j:(j,i)\in\overline{\mathcal{E}}}\langle\lambda_{(ji)},E_{ij}x_{(i)}\rangle$. The functions $q_i(\lambda)$ depend on the dual variable $\lambda$ only and the infimum is taken over local variables $x_{(i)}$. The set of minimizers of $l_i(x_{(i)},\lambda)$ with respect to $x_{(i)}$ is denoted by $X_{i,\lambda}\subseteq\mathbb{R}^{n_i}$. This allows the agents to locally compute $q_i(\lambda)$. The dual problem is given by
\begin{equation}\label{DP}
\underset{\lambda\in\mathbb{R}^{\bar{n}}}{\text{max}}\ Q(\lambda).
\end{equation}
The set of maximizers of $Q(\lambda)$ is denoted by $\Lambda_{\mathrm{opt}}$. The problem (\ref{DP}) can be solved in a standard distributed and asynchronous manner\cite{BV}. We make the following assumption.% throughout the paper.
%Each agent individually solves problem (\ref{dualind}) with respect to the primal variable $x_{(i)}$ up to certain accuracy and exchanges its solution to its neighbors asynchronously to perform the dual updates. The above process is repeated until the dual variable converges.
\begin{assumption}\label{a1}
	The following hold for Problem (\ref{PP}):
	\begin{enumerate}[(i)]
		\item Each individual cost function $F_i$ is convex.
		\item The constraint set $\mathcal{X}$ is such that (\ref{PP}) is strictly feasible.
		\item The set $X_{i,\lambda}\subseteq\mathbb{R}^{n_i}$ is non-empty for any $\lambda\in\mathbb{R}^{\bar{n}}$.
		\item The function $Q(\lambda)$ is radially unbounded: $Q(\lambda)\rightarrow-\infty$ as $|\lambda|\rightarrow\infty$. 	\hfill $\blacksquare$
	\end{enumerate} 
\end{assumption}
The convexity of functions $F_i$ in item (i) ensures that $F$ is also convex. Thus, problem (\ref{PP}) is convex with linear equality constraints. This together with item (ii) of Assumption~\ref{a1} ensures Slater's conditions hold \cite{BV}. As a result, $F^*=Q^*:=Q(\lambda^*)=\underset{\lambda}{\text{max}}Q(\lambda)$, for $\lambda^*\in\Lambda_{\mathrm{opt}}$. Item (iii) of Assumption~\ref{a1} ensures that each agent can solve its individual problem (\ref{dualind}) and return a finite solution. For example, this is the case if $\mathcal{X}\subset\mathbb{R}^n$ is chosen to be compact, see \cite{TTM} and \cite{NCN}. Item (iv) of Assumption~\ref{a1} requires $Q$ to be radially unbounded ensuring that $\Lambda_{\mathrm{opt}}$ is a compact set. The compactness will play a key role in proving Theorems~\ref{mainresult1D} and~\ref{mainresultND}. See Remark~\ref{round2r2} for a discussion on relaxing this assumption. 

\begin{remark}\label{rIman}
Differentiability of (\ref{PP}) does not play a role in our convergence analysis. This is common to many dual based algorithms analyses (see \cite{TTM}, \cite{NCN} for example) where the convergence results are provided without taking the intermediate primal minimization steps. Moreover, the proof technique used in this work does not rely on the differentiability of the dual function. This enables us to show the convergence of the algorithm in the absence of strict convexity of $F(x)$ and compactness of $\mathcal{X}$ which are typically used to guarantee differentiability of the dual function in \cite{TTM}.
	\hfill $\blacksquare$
\end{remark}
\begin{remark}\label{revise2}
% 	The main theorem establishes the convergence of the dual variable in the absence of strict convexity of $F(x)$ and differentiability of the dual function where the updates are carried out asynchronously. If the dual problem is solved, then Assumption~\ref{a1} ensures that the optimal cost $F^*=Q^*$ is found. 
	From a practical point of view, in addition to convergence of the dual variable and the optimal cost $F^*$, one may be also interested in a solution to (\ref{PP}) that is primal feasible. If this is the case, extra care must be taken to recover optimal primal solutions after the dual problem is solved. This is because for a given $\lambda\in\Lambda_{\mathrm{opt}}$, solving (\ref{dualind}) may generate $x$ that is primal-infeasible. To solve this problem, certain averaging techniques can be used as done in \cite{NO2}, \cite{FMGP} and \cite{SJR}. The simplest way is to add a quadratic regularization term with a small coefficient so that the cost function is strictly convex and the regularization term has negligible effects on the solution to the original optimization problem\cite{XJB}. %Consider a simple example of two agents with $F_1(x)=F_2(x)=\max\{|x|-1,0\}$. In general, after the dual problem is solved, it is possible to that the two agents obtain different $x$ in $[-1,1]$. To solve this, one can add $\frac{1}{1000}|x|^2$ to both $F_1$ and $F_2$, then it can be seen that both agents will obtain the solution $x_{(1)}=x_{(2)}=0$ as desired.
	\hfill $\blacksquare$
\end{remark}
\begin{remark}\label{revise3}
	The last item of Assumption~\ref{a1} is needed in proving convergence of the dual variable for our algorithm. It is commonly used in Lyapunov based proofs to conclude \textit{global} stability of a dynamical system as it ensures the existence of a valid Lyapunov function in the proof of Theorem~\ref{mainresult1D} (see Remark~\ref{round2r2} for more discussion). 
	\hfill $\blacksquare$
\end{remark}
In dual algorithms, often the explicit form of the dual function is not available. As a result, checking (iv) of Assumption~\ref{a1} \textit{a priori} may be hard. Nevertheless, sufficient conditions for radial unboundedness of the dual function as described in the following lemma can be stated 
%Before we state the lemma, we define the following set for $i\neq j$, $\bar{\mathcal{X}}_{ij}:=\{x\in\mathbb{R}^{\bar{n}_{ij}}|x=E_{ij}y,y\in\mathcal{X}_i\}$.
\begin{lemma}
	 Given $\lambda$, denote the minimizer of $l_i(x_{(i)},\lambda)$ in $\mathcal{X}_i$ by $x_{\lambda(i)}^*$. Let for all $i$, $g_i(x_{\lambda(i)}^*)\in\partial F_i (x_{\lambda(i)}^*)$ such that $g_i(x_{\lambda(i)}^*)+\sum_{j:(i,j)\in\overline{\mathcal{E}}}E_{ij}^T\lambda_{(ij)}-\sum_{j:(j,i)\in\overline{\mathcal{E}}}E_{ij}^T\lambda_{(ji)}=0$. It is further assumed that $\forall (i,j)\in\mathcal{E}$ the constraints $E_{ij}x_{(i)}=E_{ji}x_{(j)}$ are linearly independent, and there exist a $\epsilon>0$ and a non-empty subset of $\bar{\mathcal{X}}_{ij}\bigcap\bar{\mathcal{X}}_{ji}$ denoted by $\mathcal{M}_{ij}$ such that $\mathcal{M}_{ij}+\epsilon\mathbb{B}\subset\bar{\mathcal{X}}_{ij}\bigcap\bar{\mathcal{X}}_{ji}$ with $\bar{\mathcal{X}}_{ij}:=\{x\in\mathbb{R}^{\bar{n}_{ij}}|x=E_{ij}y,y\in\mathcal{X}_i\}$ and $\mathbb{B}\subset\mathbb{R}^{\bar{n}_{ij}}$. Then, $Q$ is radially unbounded. 
	\hfill $\blacksquare$
\end{lemma}
\begin{proof}
	 For any $\bar{\lambda}$, 
	\begin{align}
	&Q(\lambda)-Q(\bar{\lambda})= \sum_{i=1}^{N}F_i(x_{\lambda(i)}^*)-\sum_{i=1}^{N}F_i(x_{\bar{\lambda}(i)}^*) \notag\\
	&+\sum_{(i,j)\in\overline{\mathcal{E}}}\langle\lambda_{(ij)},E_{ij}x_{\lambda(i)}^*-E_{ji}x_{\lambda(j)}^*\rangle \notag\\
	&-\sum_{(i,j)\in\overline{\mathcal{E}}}\langle\bar{\lambda}_{(ij)},E_{ij}x_{\bar{\lambda}(i)}^*-E_{ji}x_{\bar{\lambda}(j)}^*\rangle \notag\\
	&\quad\quad\leq\sum_{i=1}^{N}g_i^T(x_{\lambda(i)}^*)(x_{\lambda(i)}^*-x_{\bar{\lambda}(i)}^*) \notag\\
	&+\sum_{j:(i,j)\in\overline{\mathcal{E}}}\langle\lambda_{(ij)},E_{ij}x_{\lambda(i)}^*\rangle-\sum_{j:(j,i)\in\overline{\mathcal{E}}}\langle\lambda_{(ji)},E_{ij}x_{\lambda(j)}^*\rangle \notag\\
	&-(\sum_{j:(i,j)\in\overline{\mathcal{E}}}\langle\bar{\lambda}_{(ij)},E_{ij}x_{\bar{\lambda}(i)}^*\rangle-\sum_{j:(j,i)\in\overline{\mathcal{E}}}\langle\bar{\lambda}_{(ji)},E_{ij}x_{\bar{\lambda}(j)}^*\rangle)\notag \\
	&\leq -\sum_{i=1}^{N}g_i^T(x_{\lambda(i)}^*)x_{\bar{\lambda}(i)}^*-C \notag\\
	&=-\sum_{(i,j)\in\overline{\mathcal{E}}}\langle\lambda_{(ij)},E_{ij}x_{\bar{\lambda}(i)}^*-E_{ji}x_{\bar{\lambda}(j)}^*\rangle-C, \notag
	\end{align}
	where the first inequality is due to the definition of subgradients, $C := \sum_{j:(i,j)\in\overline{\mathcal{E}}}\langle\bar{\lambda}_{(ij)},E_{ij}x_{\bar{\lambda}(i)}^*\rangle-\sum_{j:(j,i)\in\overline{\mathcal{E}}}\langle\bar{\lambda}_{(ji)},E_{ij}x_{\bar{\lambda}(j)}^*\rangle$, and the last inequality comes from the assumption $g_i(x_{\lambda(i)}^*)+\sum_{j:(i,j)\in\overline{\mathcal{E}}}E_{ij}^T\lambda_{(ij)}-\sum_{j:(j,i)\in\overline{\mathcal{E}}}E_{ij}^T\lambda_{(ji)}=0$. Due to the assumption that $\mathcal{M}_{ij}+\epsilon\mathbb{B}\subset\bar{\mathcal{X}}_{ij}\bigcap\bar{\mathcal{X}}_{ji}$, the case where a component of $\lambda$ is identically multiplied by $0$ does not exist. It also ensures that there are elements in $\mathcal{X}$ such that the signs for all components of $\lambda$ can be made negative or positive as one wishes by selecting appropriate $x\in\mathcal{X}$. The question now becomes for such $x\in\mathcal{X}$, can we find the corresponding $\bar{\lambda}$. It is obvious that, all linearly independent constraints can be put in a compact form of $Hx=0$. Therefore, if for any $x$, we can find a corresponding $\bar{\lambda}$ such that $g(x)=-H^T\lambda$, where $g(x)\in\partial F(x)$, then the desired $x_{\bar{\lambda}(i)}^*$ always exists for any $i$. Indeed, since we assume that the constraints are linearly independent, $H$ has full row rank and $H^T$ has full column rank. Therefore, by choosing appropriate $\bar{\lambda}$, $-\sum_{(i,j)\in\overline{\mathcal{E}}}\langle\lambda_{(ij)},E_{ij}x_{\bar{\lambda}(i)}^*-E_{ji}x_{\bar{\lambda}(j)}^*\rangle$ tends to $-\infty$ as $|\lambda|$ goes to infinity. Thus $Q$ is radially unbounded.
\end{proof}
\section{The Asynchronous Block Coordinate Subgradient Algorithm}\label{sec3}
In this section, we present an asynchronous block supergradient algorithm to solve the dual problem (\ref{DP}). First, let us re-write the Lagrangian (\ref{lagrangian}) in a more compact form:
\begin{equation}\label{compactlagrangian}
L(x, \lambda):=F(x)+\langle\lambda,\mathbf{E}(x)\rangle,
\end{equation}
where the linear constraint $\mathbf{E}(x):=(\ldots,\mathbf{E}_i(x),\ldots)$ is given by the concatenation of matrices $\mathbf{E}_i(x):=(\ldots,E_{ij}x_{(i)}-E_{ji}x_{(j)},\ldots),(i,j)\in\overline{\mathcal{E}}$ such that $i\in\mathcal{I}$. For $\lambda\in\mathbb{R}^{\bar{n}}$, let $x_\lambda^*=(x_{\lambda(1)}^*,x_{\lambda(2)}^*,\ldots,x_{\lambda(N)}^*)$ be an arbitrary real minimizer of the Lagrangian $L(x, \lambda)$, where $x_{\lambda(i)}^*\in\mathcal{X}_{i},i\in\mathcal{V}$ are minimizers of $l_i(x_{(i)},\lambda)$. We next state one way to compute a supergradient of $Q(\lambda)$. The proof can be found in standard textbooks, e.g., \cite{Bert1}, and is omitted.
% We next give the result showing that a supergradient of $Q(\lambda)$ can be obtained by minimizing $L(x, \lambda)$ with respect to $x$. The proof can be found in standard textbooks, e.g., \cite{Bert1}, and is omitted.
\begin{lemma}\label{superg}
	Let $Q(\lambda)$ be the dual function (\ref{dual}), then for any $\lambda\in\mathbb{R}^{\bar{n}}$ and a corresponding $x^*_{\lambda}\in\mathcal{X}$, $\mathbf{E}(x_\lambda^*)\in\partial^+Q(\lambda)$. \hfill $\blacksquare$
\end{lemma}
%\begin{proof}
%	For $\lambda,\bar{\lambda}\in\mathbb{R}^{\bar{n}}$, we have: $Q(\bar{\lambda})=\underset{x\in\mathbb{R}^n}{\text{inf}}L(x, \bar{\lambda})=\underset{x\in\mathbb{R}^n}{\text{inf}}F(x)+\langle\bar{\lambda},\mathbf{E}(x)\rangle\leq F(x_\lambda^*)+\langle\bar{\lambda},\mathbf{E}(x_\lambda^*)\rangle=F(x_\lambda^*)+\langle\lambda,\mathbf{E}(x_\lambda^*)\rangle+\langle\bar{\lambda}-\lambda,\mathbf{E}(x_\lambda^*)\rangle=Q(\lambda)+\langle\bar{\lambda}-\lambda,\mathbf{E}(x_\lambda^*)\rangle$. By definition, $\mathbf{E}(x_\lambda^*)$ belongs to $\partial Q(\lambda)$.
%\end{proof}

This implies $E_{ij}x_{\lambda(i)}^*-E_{ji}x_{\lambda(j)}^*$ is a component of the supergradient associated with the corresponding multiplier. The synchronous supergradient update for $\lambda$ is given by
\begin{equation}\label{dualasc}
\lambda^+=\lambda+\alpha g(\lambda)=\lambda+\alpha\mathbf{E}(x_\lambda^*),
\end{equation}
where $\alpha>0$ is the stepsize to be appropriately chosen and $g\in\partial^+ Q$. As a result, asynchronous version of the dual updates corresponds to the supergradient update (\ref{dualasc}) for some of its components
$\lambda_{(ij)}^+=\lambda_{(ij)}+\alpha(E_{ij}x_{\lambda(i)}^*-E_{ji}x_{\lambda(j)}^*)$,
for some $(i,j)\in\overline{\mathcal{E}}$. Specifically, each agent minimizes $l_i(x_{(i)},\lambda)$ with respect to $x_{(i)}$ locally. When a pair of agents finish solving their local problems, they can communicate their solutions to each other to perform a dual update without waiting for others to complete their computations. These communication times may depend on agents' computational capabilities, availability of the network and so on. Based on Lemma~\ref{superg}, the proposed dual asynchronous algorithm that generates the sequence $\{\lambda_{[k]}\}_{k=1}^{\infty}$ and $\{\alpha_{[k]}\}_{k=1}^{\infty}$ can be modelled by the following stochastic discrete-time system
\begin{subequations}
\label{reducedmodelfull}
\begin{align}
\lambda^{+}&\in \lambda+\Phi(\alpha\odot v^+)\odot(\partial^+ Q(\lambda)+e^+)\\
\alpha^{+}&=A(\alpha)\odot v^++\alpha\odot(\mathbf{1}_{|\overline{\mathcal{E}}|}-v^+)\\
\gamma^{+}&=\gamma+v^+,
\end{align}
\end{subequations}
where $\odot$ denotes the Hadamard product of two matrices, $\alpha=(\ldots,\alpha_{(ij)},\ldots)\in\mathbb{R}_{>0}^{|\overline{\mathcal{E}}|}$, $(i,j)\in\overline{\mathcal{E}}$ is the stepsize vector, $\gamma=(\ldots,\gamma_{(ij)},\ldots)\in\mathbb{R}_{\geq0}^{|\overline{\mathcal{E}}|}$, $(i,j)\in\overline{\mathcal{E}}$ is the vector of the number of dual updates up to the current time instance and its initial condition is therefore set to $0$. The function $A(\alpha)=(\ldots,A_{(ij)}(\alpha_{(ij)}),\ldots)\in\mathbb{R}_{>0}^{|\overline{\mathcal{E}}|}$, $(i,j)\in\overline{\mathcal{E}}$ is the vector representing the dynamics of local stepsize update. To model the asynchrony, we introduce the random input $v=(\ldots,v_{(ij)},\ldots)$ such that $v_{(ij),{[k+1]}}$ is equal to $1$ if at time $k$ edge $(i,j)\in\overline{\mathcal{E}}$ is updating and is equal to $0$ otherwise. Furthermore, $\Phi(\alpha\odot v^+):=(\ldots,\alpha_{(ij)}v^+_{(ij)}\mathbf{1}_{\bar{n}_{ij}},\ldots)$. The other random input $e^+=(\ldots,e_{(ij)},\ldots)$ following some probability distribution satisfying Assumption~\ref{a6} accounts for supergradient evaluation errors coming from communication noise and the fact that the local problems may not be solved exactly in finite time. For notational convenience, we re-label the vectors $\lambda$, $\alpha$ and $\gamma$ in (\ref{reducedmodelfull}) and give a component-wise representation of (\ref{reducedmodelfull}). Note that there are $|\overline{\mathcal{E}}|$ dual updates and we label the corresponding update variables from $1$ to $|\overline{\mathcal{E}}|$, i.e., $\lambda=(\lambda_{(1)},\ldots,\lambda_{(|\overline{\mathcal{E}}|)})$, $\alpha=(\alpha_{(1)},\ldots,\alpha_{(|\overline{\mathcal{E}}|)})$ and $\gamma=(\gamma_{(1)},\ldots,\gamma_{(|\overline{\mathcal{E}}|)})$. The mappings $A(\cdot)$ and $\partial^+ Q(\cdot)$ are partitioned similarly. For $i\in\mathbb{Z}_{|\overline{\mathcal{E}}|}\backslash\{0\}$, the $i^{th}$ component of (\ref{reducedmodelfull}) can be written as 
\begin{subequations}\label{rmcomponent}
\begin{align}
\lambda_{(i)}^{+}&\in\lambda_{(i)}+\alpha_{(i)}v^+_{(i)}(\partial^+ Q(\lambda)_{(i)}+e_{(i)}^+)\\
\alpha_{(i)}^{+}&=A_{(i)}(\alpha_{(i)})v^+_{(i)}+\alpha_{(i)}(1-v^+_{(i)})\\
\gamma_{(i)}^+&=\gamma_{(i)}+v^+_{(i)}.
\end{align}
\end{subequations}

An implementation of \eqref{reducedmodelfull} is presented in Algorithm~\ref{alg1}.
\begin{algorithm} 
	
	\caption{Asynchronous Block Supergradient Algorithm}\label{alg1}
	\begin{algorithmic}[1]
		\State \textbf{Initialization}: $\lambda_{(ij),[0]}$ $\forall (i,j)\in\overline{\mathcal{E}}$.
		\State \textbf{Primal Update}: At time $k$, for all agents $i\in\mathcal{V}$,
		\begin{equation*}
		x_{(i)}^*\leftarrow \; \underset{{x_{(i)}\in \mathcal{X}_i}}{\arg\min}\ l_i(x_{(i)},\lambda)
		\end{equation*}
		\State \textbf{Dual Update}: For all $(i,j)\in\overline{\mathcal{E}}$ such that $v_{(ij),[k+1]}=1$:
		\begin{align*}
		\lambda_{(ij)}^+&\leftarrow\lambda_{(ij)}+\alpha_{(ij)}(E_{ij}x_{(i)}^*-E_{ji}x_{(j)}^*+e_{(ij)}^+),\\
		\alpha_{(ij)}^+&\leftarrow A_{ij}(\alpha_{(ij)}).
		\end{align*}
		For all $(i,j)\in\overline{\mathcal{E}}$ such that $v_{(ij),[k+1]}=0$:
		\begin{equation*}
		\lambda_{(ij)}^+\leftarrow\lambda_{(ij)},\quad 
		\alpha_{(ij)}^+\leftarrow \alpha_{(ij)}.
		\end{equation*}
		\State Set $k\leftarrow k+1$ and go to Step 2.
	\end{algorithmic}
\end{algorithm}

We make the following assumption on the update protocol.
\begin{assumption}\label{a4}
	 For any $i\in\mathbb{Z}_{|\overline{\mathcal{E}}|}\backslash\{0\}$ and $m\in\mathbb{Z}_{\geq 0}$, there exists a $\Delta>0$ such that $\liminf_{k\rightarrow\infty}\frac{{\gamma}_{(i),[k+m]}}{k+m}\geq\Delta$ almost surely, where ${\gamma}_{(i)}$ is generated by (\ref{rmcomponent}) with ${\gamma}_{(i),[0]}=0$.
	\hfill $\blacksquare$
\end{assumption}
\begin{remark}\label{rdis}
	Assumption~\ref{a4} requires that each coordinate is updated infinitely often. As shown in \cite{Borkar98}, it can be satisfied whenever each coordinate, given the past history, has an updating probability lower bounded by a positive constant\cite{WO}. Other examples of stochastic updating conditions satisfying Assumption~\ref{a4} include the case where coordinates are updated in an i.i.d.~fashion \cite{LSN}, or when the coordinates are updated based on a controlled Markov chain \cite{Borkar1}. %such that any stationary policy results in an irreducible chain with a stationary distribution assigning probabilities lower bounded by a positive constant to each coordinate \cite{Borkar1}. 
	Two examples for deterministic rules satisfying Assumption~\ref{a4} include cyclic updating of coordinates \cite{Wright} and the case of uniformly persistently exciting updates \cite{TN} where each coordinate is updated at least once in $N$ iterations for a given integer $N$ \cite{FGGN}. 
	\hfill $\blacksquare$
\end{remark}
For $x\in\mathbb{R}$, let $\lfloor x\rfloor$ denote the maximum integer no larger than $x$. We present next the assumptions on agent's local stepsize. 
\begin{assumption}\label{a5}
	For any $i,j\in\mathbb{Z}_{|\overline{\mathcal{E}}|}\backslash\{0\}$ and stepsize sequences $\{\alpha_{(i)k}\}$ and $\{\alpha_{{(j)}k}\}$ generated by $\alpha_{(i)}^+=A_{(i)}(\alpha_{(i)})$ and $\alpha_{(j)}^+=A_{(j)}(\alpha_{(j)})$ respectively, there exist positive scalars $\bar{\delta}$ and $\bar{\Delta}>0$ such that {\it(i)} $\sum_{k=0}^{\infty}\alpha_{(i),[k]}=\infty$ and $\sum_{k=0}^{\infty}\alpha_{(i),[k]}^2<\infty$; {\it(ii)} $\limsup_{k\rightarrow\infty}\sup_{\delta_1\in[\delta_2,1]}\alpha_{{(i)},[{\lfloor\delta_1 k\rfloor}]}/\alpha_{{(j)},[k]}\leq\bar{\Delta}$, for some $\delta_2 \in (0,\Delta]$;
		{\it(iii)} $\liminf_{k\rightarrow\infty}\inf_{\delta_1\in[\delta_2,1]}\alpha_{{(i)},[{\lfloor\delta_1 k\rfloor}]}/\alpha_{{(j)},[k]}\geq\bar{\delta}$, for some $\delta_2 \in (0,\Delta]$.
		\hfill $\blacksquare$
\end{assumption}
\begin{remark}\label{r4}
	Since we are not working exclusively with Lipschitz continuous gradients, it is reasonable to use diminishing stepsizes to achieve almost sure convergence to the exact solution set $\Lambda_{\mathrm{opt}}$. Item (ii) of Assumption~\ref{a5} requires the stepsizes of agents in the network to be roughly in the same ``time scale''. Moreover, different from \cite{LSN}, we allow agents to use different local stepsizes whose updates are not required to be coordinated as long as they satisfy Assumption~\ref{a5}. They can be easily satisfied by choosing commonly used stepsizes including $\{c_i/k\}$ and $\{c_i/(k\log k)\}$ (with appropriate modifications when $k=0$) where $c_i>0$ is a bounded constant that agent $i\in\mathcal{I}$ can choose locally. 
	\hfill $\blacksquare$
\end{remark}
%This is to avoid situations where the updates on one particular coordinate goes too quick. For example, suppose two coordinates of the dual function are updated respectively using stepsizes $\{1/n\}$ and $\{1/n^\frac{3}{4}\}$. It is clear that the stepsizes $\{1/n^\frac{3}{4}\}$ will lead to a faster ``time scale'' meaning the progress on the other coordinate is negligible. The analysis of the convergence properties in this situation may be done via analysing stability properties of the corresponding limiting continuous time system which is singularly perturbed \cite{KKO} due to the time scale separation induced by the stepsizes. However, this is beyond the scope of this work and we refer readers to Chapter 6 of \cite{Borkar1} for more details.
We denote the history of (\ref{reducedmodelfull}) and (\ref{rmcomponent}) up to time $k\in\mathbb{Z}_{\geq0}$ by $\mathcal{F}_k:=\sigma(\lambda_{[m]},\alpha_{[m]},e_{[m]},v_{[m]},m\in\mathbb{Z}_k)$, which is a $\sigma$-subfield of $\mathcal{F}$. The following assumption on the supergradient error states that the error is a martingale difference noise as shown later in the proof of Theorem~\ref{mainresult1D} which is a standard assumption in the stochastic approximation literature \cite{Borkar1}.
\begin{assumption}\label{a6}
	There exist a sequence $\{\beta_{[k]}\}_{k=0}^{\infty}$ and a $K\geq0$ such that the random sequence of error vectors $e$ in (\ref{reducedmodelfull}) satisfies $\mathbb{E}[e_{[k+1]}|\mathcal{F}_k]=\beta_{[k]}$ and $\mathbb{E}[|e_{[k+1]}|^2|\mathcal{F}_k]\leq K$.
	\hfill $\blacksquare$
\end{assumption}

\begin{remark}\label{revised4}
	Most stochastic approximation and optimization literature considers zero mean martingale difference noise under which exact almost sure convergence can be proved. However, in some applications, it is natural that the noise may not be zero mean. We introduce the variable $\beta_{[k]}$ to cover both the zero mean ($\beta_{[k]}=0$) and non-zero mean ($\beta_{[k]}\neq0$) cases. We prove that if $\beta_{[k]}=0$ then we can establish the same convergence results as shown in the stochastic approximation literature. If $\beta_{[k]}\neq0$, then we can show convergence to a neighborhood of $\Lambda_\mathrm{opt}$.
	\hfill $\blacksquare$
\end{remark}
\begin{assumption}\label{arevised1}
	The dual variable $\lambda_{[k]}$ is bounded almost surely for any $k$.
	\hfill $\blacksquare$
\end{assumption}
Assumption~\ref{arevised1} is standard in stochastic approximation literature \cite{Borkar1} to guarantee almost sure convergence. In the next section, we will give a sufficient condition for satisfaction of Assumption~\ref{arevised1}. 

The main result of the paper is stated below.
\begin{theorem}\label{mainresult1D}
	Let Assumptions~\ref{a1}--\ref{arevised1} hold. If $\beta_{[k]}=0$ for all $k\in\mathbb{Z}_{\geq0}$, then $\lim_{k\rightarrow\infty}|\lambda_{[k]}|_{\Lambda_{\mathrm{opt}}}=0$ almost surely for any $(\lambda_{[0]},\alpha_{[0]},\gamma_{[0]})\in\mathbb{R}^{\bar{n}}\times\mathbb{R}_{>0}^{|\overline{\mathcal{E}}|}\times\{0\}$. 	\hfill $\blacksquare$
\end{theorem}
The proof of Theorem~\ref{mainresult1D} is inspired by the asynchronous stochastic approximation idea introduced in \cite[Chapter 7]{Borkar1} and \cite{Borkar98} where similar results are established for the case where the cost function is differentiable. Before presenting it in detail, a brief outline of the proof followed by necessary intermediate results are provided below.
% The proofs of intermediate results can be found in the Appendix. 

The proof generalises that of \cite[Chapter 7, Theorem 2]{Borkar1} to the non-smooth case. We first linearly interpolate the asynchronous dual updates to make it a ``continuous time'' signal. We then show that this continuous time signal almost surely converges to a solution of a non-autonomous differential inclusion (Proposition~\ref{Prop2}) with $\Lambda_{\mathrm{opt}}$ as its uniformly globally asymptotically stable (UGAS) positively invariant set in the sense of \cite[Definition 2.2]{LSW} (Proposition~\ref{limi}). %The necessary preliminaries and proofs of the intermediate results will be given in the Appendix.
% \begin{proof}
	 For notational convenience, we define the radially unbounded convex function $f:=-Q$ whose minimisers form a non-empty compact set (Assumption~\ref{a1}). Hence, for $i\in\mathcal{I}$, (\ref{rmcomponent}) can be written as
	\begin{subequations}\label{e1.1a}
	\begin{align}
	\lambda_{(i)}^{+}&= \lambda_{(i)}-\alpha_{(i)}v^+_{(i)}g(\lambda)_{(i)}+\alpha_{(i)}v^+_{(i)}e_{(i)}^+\\
	\alpha_{(i)}^{+}&=A_{(i)}(\alpha_{(i)})v^+_{(i)}+\alpha_{(i)}(1-v^+_{(i)})\\
	\gamma_{(i)}^+&=\gamma_{(i)}+v^+_{(i)},
	\end{align}
	\end{subequations}
	where $g(\lambda)\in\partial f(\lambda)$, and is calculated by taking the difference of the two corresponding primal solutions of local problems and can be any element of $\partial f(\lambda)$. 
	
	Note that for each $i$, $\lambda_{(i)}$ in (\ref{e1.1a}) can be written as
	\begin{equation}\label{e2.3}
	\lambda_{(i)}^+=\lambda_{(i)}-\bar{\alpha}\tau_{(i)}(g(\lambda)_{(i)}-e_{(i)}^+),
	\end{equation}
	where $\bar{\alpha}:=\underset{i\in \mathcal{I}_a}{\max}\ \alpha_{(i)}$, $\tau_{(i)}:=\frac{\alpha_{(i)}}{\bar{\alpha}}v^+_{(i)}$ and $\mathcal{I}_a\subseteq\mathbb{Z}_{|\overline{\mathcal{E}}|}\backslash\{0\}$ is the set of components that are updating (we define $\bar{\alpha}$ and $\tau_{(i)}$ to be $0$ when $\mathcal{I}_a=\emptyset$). The next lemma shows that the sequence $\{\bar{{\alpha}}_{[k]}\}$ is square summable but not summable.
	\begin{lemma}\label{newstepsize}
		The sequence $\{\bar{{\alpha}}_{[k]}\}$ satisfies $\sum_{k=0}^{\infty}\bar{{\alpha}}_{[k]}^2<\infty$ and $\sum_{k=0}^{\infty}\bar{{\alpha}}_{[k]}=\infty$ almost surely. \hfill $\blacksquare$
	\end{lemma}
  Let 
  \begin{align}\label{eq:t_k}
     t_{[k]}:=\sum_{m=0}^{k-1}\bar{{\alpha}}_{[m]}, \ k\geq 1, \ t_{[0]}=0, 
  \end{align} and for $t\in[t_{[k]},t_{[k+1]})$% by setting $\bar{\lambda}(t_{[k]})=\lambda_{[k]}$ and
  \begin{equation}\label{e2.4}
  \bar{\lambda}(t):=\lambda_{[k]}+\frac{\lambda_{[k+1]}-\lambda_{[k]}}{\bar{{\alpha}}_{[k]}}(t-t_{[k]}), \ \bar{\lambda}(t_{[k]})=\lambda_{[k]}.
  \end{equation} 
  % the linearly interpolated process $\bar{\lambda}(t)$ on the interval $[t_{[k]},t_{[k+1]})$
  By construction, the linearly interpolated process $\bar{\lambda}(t)$ is continuous and piecewise linear. Let $\theta(t):=\text{diag}(\tau_{(1),[k]},\ldots,\tau_{(|\overline{\mathcal{E}}|),[k]}), \text{for}\ t\in[t_{[k]},t_{[k+1]})$. Intuitively, the sum of $\tau_{(i)}$ indicates how frequently the agent $i$ is updating. By Lemma~\ref{LL} (see Appendix), $f(\lambda)$ is guaranteed to be locally Lipschitz. Thus by Rademacher's theorem \cite{Clarke}, $f(\lambda)$ is differentiable almost everywhere. Hence, for almost all $\lambda\in\mathbb{R}^{\bar{n}}$, $\nabla f(\lambda)$ is well-defined, continuous, and $\partial f(\lambda)=\{\nabla f(\lambda)\}$. Therefore, $g(\lambda)$ is continuous almost everywhere. Equivalently, by defining $r(t):=\max\{k\in\mathbb{Z}_{\geq0}:t_{[k]}\leq t\}$: %the following relationship on the interpolated process
  \begin{align}
  \bar{\lambda}(t_{[k]}+u)-\bar{\lambda}(t_{[k]})=&-\int_{0}^{u}\theta(t_{[k]}+s)g(\lambda_{r(t_{[k]}+s)})ds \notag \\
  &+{\zeta}_{r(t_{[k]}+u)}-{\zeta}_{r(t_{[k]})}, \label{e2.5}
  \end{align}
  where ${\zeta}_{[r(t_{[k]})]}={\zeta}_{[k]}:=\sum_{m=0}^{k-1}\bar{{\alpha}}_{[m]}\theta(t_{[m]})e_{[m+1]}, k\geq 1$. Compared to $\bar{\lambda}(t)$, which is piecewise linear and continuous, $\lambda_{[r(t)]}$ is piecewise constant and may be discontinuous at $t=t_{[k]}$. From Assumption~\ref{a6}, it can be seen that $({\zeta}_{[k]},\mathcal{F}_k)$ is a martingale process. Moreover, the following proposition states that if $\beta_{[k]}=0$ for all $k$, the effect of the term ${\zeta}_{r(t_{[k]}+u)}-{\zeta}_{r(t_{[k]})}$ in (\ref{e2.5}) diminishes as $k\rightarrow\infty$.
  \begin{proposition}\label{martingaleconvergence}
  	If $\beta_{[k]}=0$ for all $k$, then $\lim_{k\rightarrow\infty}\zeta_{[k]}<\infty$ and the sequence of functions of $u$, $\{{\zeta}_{[r(t_{[k]}+u)]}-{\zeta}_{[r(t_{[k]})]}\}_{k=1}^{\infty}$ converges to $0$ on $[0,\infty)$ almost surely. \hfill $\blacksquare$
  \end{proposition}
  Next, by defining $\phi(t+u,t,\bar{\lambda}(t)):=\bar{\lambda}(t)-\int_{0}^{u}\theta(t+s)g(\bar{\lambda}(t+s))ds$, we have the following result which shows $\phi(t+u,t,\bar{\lambda}(t))$ can be used to approximate $\bar{\lambda}(t+u)$ when $t$ is sufficiently large.
  \begin{proposition}\label{lemma3}
  	For any $T>0$, $\lim_{t\rightarrow\infty}\sup_{u\in[0,T]}|\bar{\lambda}(t+u)-\phi(t+u,t,\bar{\lambda}(t))|=0$ almost surely. \hfill $\blacksquare$
  \end{proposition}
  The behavior of the limit point of $\{\bar{\lambda}(t_{[k]}+\cdot)\}$ is studied next. 
  \begin{proposition}\label{Prop2}
  	 For any positive scalar $T$ and positive integer $\bar{n}$, any limit point $\tilde{\lambda}(\cdot)$ of $\{\bar{\lambda}(t_{[k]}+\cdot)\}$ defined in (\ref{e2.4}) in the space of continuous functions from $[0,T]$ to $\mathbb{R}^n$ denoted by $C([0,T];\mathbb{R}^{\bar{n}})$ is almost surely a solution to the non-autonomous differential inclusion 
  	\begin{equation}\label{e2.1}
  	\tilde{\lambda}\in-\Theta(t)\partial f(\tilde{\lambda}),
  	\end{equation}
  	where $\Theta(t)$ is a diagonal matrix-valued measurable function with diagonal entries in $[\varepsilon,1]$, with $0<\varepsilon\leq1$. 
  	\hfill $\blacksquare$
  \end{proposition}
  Lastly, the following result states stability properties of the limiting differential inclusion (\ref{e2.1}) which also ensures that no finite escape time exists for (\ref{e2.1}) by compactness of $\Lambda_{\mathrm{opt}}$.
  \begin{proposition}\label{limi}
  	Suppose $\Theta(t)$ is a diagonal matrix-valued measurable function with diagonal entries in $[\varepsilon,1]$, then (\ref{e2.1}) is UGAS with respect to the compact set $\Lambda_{\mathrm{opt}}$.
  	\hfill $\blacksquare$
  \end{proposition}
  Now, we are ready to prove Theorem~\ref{mainresult1D}.
  % \todo[inline]{Your proof need to be a sequence of $P\implies Q$ while invoking the necessary Lemma/Proposition along the way. This has to be rewritten.}

\noindent{\it Proof of Theorem~\ref{mainresult1D}:}
  By Lemma~\ref{newstepsize}, $t_{[k]}\rightarrow\infty$ (defined in \eqref{eq:t_k}) as $k\rightarrow\infty$. Thus, $\bar{\lambda}(t)$ is defined over all $t\in[0,\infty)$. From \eqref{e2.4}, the limit $\lambda_{[k]}$ as $k \rightarrow \infty$ exists if and only if the limit of $\bar{\lambda}(t)$ as $t\rightarrow \infty$ exists. These limits are the same.
  % ?{ why thus?}, when $k\rightarrow\infty$, the behavior of $\lambda_{[k]}$ can be characterized { what does can be characterise mean? I think all of these should be moved to the end where you invoke Borkar's result} by the interpolated signal $\bar{\lambda}(t)$ when $t\rightarrow\infty$ . 
  By Lemma~\ref{LL}, $\partial f(\lambda)$ is locally bounded. Moreover, since the sequence $\{\lambda_{[k]}\}_{k=0}^{\infty}$ generated by (\ref{e1.1a}) is almost surely bounded, $g(\lambda_{[k]})$ is almost surely bounded for all $g(\lambda_{[k]})\in\partial f(\lambda_{[k]})$, $k\in\mathbb{Z}_{\geq0}$. Furthermore, $\theta(t)$ is uniformly bounded for all $t\in[0,\infty)$. Combining these facts and Proposition~\ref{martingaleconvergence}, it can be seen that, almost surely, $\bar{\lambda}(t_{[k]}+u)-\bar{\lambda}(t_{[k]})$ in (\ref{e2.5}) is norm bounded. Hence, the sequence of continuous time processes $\{\bar{\lambda}(t_{[k]}+\cdot)\}_{k=0}^{\infty}$ is equicontinuous (see Preliminaries in the appendix for the definition) and bounded. Note that, for a given $k$, $\bar{\lambda}(t_{[k]}+\cdot)$ is a function defined on $[0,\infty)$ shifted by $t_{[k]}$, thus the sequence of functions $\{\bar{\lambda}(t_{[k]}+\cdot)\}_{k=0}^{\infty}$ are in fact defined on the same interval, but shifted differently. Then by the Arzel\`{a}-Ascoli Theorem (see Lemma~\ref{AA} in Appendix), for any $T>0$, $\{\bar{\lambda}(t_{[k]}+\cdot)\}_{k=0}^{\infty}$ is relatively compact in $C([0,T];\mathbb{R}^{\bar{n}})$. Thus, there exists a subsequence of $\{\bar{\lambda}(t_{[k]}+\cdot)\}_{k=0}^{\infty}$ that converges to a limit point. Following similar arguments, it can be shown that, for $t\geq0$, $\{\phi(t_{[k]}+\cdot,t_{[k]},\bar{\lambda}(t_{[k]}))\}$ is also an equicontinuous, point-wise bounded family of trajectories.
  Therefore, it follows from Proposition~\ref{lemma3} that $\{\phi(t_{[k]}+\cdot,t_{[k]},\bar{\lambda}(t_{[k]}))\}$ and $\{\bar{\lambda}(t_{[k]}+\cdot)\}$ have the same limit points in $C([0,T];\mathbb{R}^{\bar{n}})$. As a result, studying the limiting behavior of $\{\bar{\lambda}(t_{[k]}+\cdot)\}$ can be done via looking at the corresponding common subsequence of $\{\phi(t_{[k]}+\cdot,t_{[k]},\bar{\lambda}(t_{[k]}))\}$. 
  By Propositions~\ref{lemma3} and~\ref{Prop2}, almost surely, all limit points of the interpolated process (\ref{e2.4}) are solutions of (\ref{e2.1}) which is UGAS with respect to $\Lambda_{\mathrm{opt}}$ as shown in Proposition~\ref{limi}. Thus, from \cite[Chapter 5, Corollary 4]{Borkar1}, $\lambda_{[k]}$ in (\ref{e1.1a}) almost surely converges to $\Lambda_{\mathrm{opt}}$.  \hfill $\blacksquare$
  
By imposing the following additional assumption on the stepsize sequence, we are able to state almost sure convergence of $\lambda_{[k]}$ to a neighborhood of $\Lambda_{\mathrm{opt}}$ if $\beta_{[k]}$ is norm-bounded.
\begin{assumption}\label{aaditional}
	For any $i,j\in\mathbb{Z}_{|\overline{\mathcal{E}}|}\backslash\{0\}$, the stepsize sequences $\{\alpha_{(i),[k]}\}$ and $\{\alpha_{{(j)},[k]}\}$ where $\alpha_{(i)}^+=A_{(i)}(\alpha_{(i)})$ and $\alpha_{(j)}^+=A_{(j)}(\alpha_{(j)})$ satisfy the following:
	\begin{enumerate}[(i)]
	  {
	  \item $(\exists c_{ij}>0)\quad \lim_{n\rightarrow\infty}\frac{\sum_{k=0}^{n}\alpha_{(i),[k]}}{\sum_{k=0}^{n}\alpha_{{(j)},[k]}}=c_{ij}$;
	  \item $(\exists b_1\in(0,\Delta])\;(\forall b_2\in[b_1,1])\; (\exists c\in(0,1]) \quad \lim_{n\rightarrow\infty}\frac{\sum_{k=0}^{\lfloor b_2n\rfloor}\alpha_{(i),[k]}}{\sum_{k=0}^{n}\alpha_{{(i)},[k]}}=c$.}
% 		\item There exists a $c_{ij}>0$ such that $\lim_{n\rightarrow\infty}\frac{\sum_{k=0}^{n}\alpha_{(i),[k]}}{\sum_{k=0}^{n}\alpha_{{(j)},[k]}}=c_{ij}$.
% 		\item There exists some $b_1\in(0,\Delta]$ such that, for any $b_2\in[b_1,1]$, there exists some constant $c\in(0,1]$ such that $\lim_{n\rightarrow\infty}\frac{\sum_{k=0}^{\lfloor b_2n\rfloor}\alpha_{(i),[k]}}{\sum_{k=0}^{n}\alpha_{{(i)},[k]}}=c$. 
		\hfill $\blacksquare$
	\end{enumerate}
\end{assumption}
\begin{theorem}\label{mainresultND}
	Let Assumptions~\ref{a1}--\ref{aaditional} hold. Suppose that for any $\delta_e>0$ and $\Delta_e>0$ there exists $\varepsilon_e>0$, such that $|\beta_{[k]}|\leq\varepsilon_e$ and the sequence $\{\lambda_{[k]}\}_{k=1}^{\infty}$ generated by (\ref{reducedmodelfull}) is bounded almost surely and $\lambda_{[k]}\in\Delta_e\mathbb{B}$ infinitely often. Then $\lim_{k\rightarrow\infty}|\lambda_{[k]}|_{\Lambda_{\mathrm{opt}}+\delta_e\mathbb{B}}=0$ almost surely for any $(\lambda_{[0]},\alpha_{[0]},\gamma_{[0]})\in\mathbb{R}^{\bar{n}}\times\mathbb{R}_{>0}^{|\overline{\mathcal{E}}|}\times\{0\}$.	 \hfill $\blacksquare$
\end{theorem}
\begin{proof}
Note that, we can split the non-zero-mean error into two terms with one being the martingale difference noise considered in Theorem~\ref{mainresult1D} and the other one being a noise term bounded by $\varepsilon_e$. Following similar arguments in the proof of Theorem~\ref{mainresult1D}, it can be seen that the limiting differential inclusion becomes $\dot{\tilde{\lambda}}\in-\Theta(t)(\partial f(\tilde{\lambda})+\varepsilon_e\mathbb{B})$. Since there exists $c_{ij}>0$ such that $\lim_{n\rightarrow\infty}\frac{\sum_{k=0}^{n}\alpha_{{(i),[k]}}}{\sum_{k=0}^{n}\alpha_{(j),[k]}}=c_{ij}$, the stepsize sequences used by agents belong to the class of balanced stepsizes defined in \cite[Section 3]{Borkar98}. Thus, following the arguments in \cite[Theorem 3.2]{Borkar98}, $\Theta(t)=\Theta$ and, hence, the deterministic limiting differential inclusion becomes 
\begin{equation}\label{LDsim}
\dot{\tilde{\lambda}}\in-\Theta(\partial f(\tilde{\lambda})+\varepsilon_e\mathbb{B}),
\end{equation}
which is autonomous. By Lemma~\ref{LL} (see Appendix), $\partial f(\tilde{\lambda})$ is outer semicontinuous, convex and locally bounded\footnote{A set-valued mapping $M:\mathbb{R}^p\rightrightarrows\mathbb{R}^n$ is outer semicontinuous if for each $(x_i,y_i)\rightarrow(x,y)\in\mathbb{R}^p\times\mathbb{R}^n$ satisfying $y_i\in M_i(x_i)$ for all $i\in\mathbb{Z}_{\geq 0}$, $y\in M(x)$. It is locally bounded if for each bounded $K\subset\mathbb{R}^p$, $M(K):=\bigcup_{x\in K}M(x)$ is bounded.}. Thus, the continuous time system $\dot{\tilde{\lambda}}\in-\Theta\partial f(\tilde{\lambda})$ satisfies \cite[Assumption 6.5]{GST} and is well-posed. Since $\Theta(\partial f(\tilde{\lambda})+\varepsilon_e\mathbb{B})\subseteq\Theta(\partial f(\tilde{\lambda}))+\varepsilon_e\mathbb{B}$ and $\Lambda_{\mathrm{opt}}$ is a compact set by Assumption~\ref{a1}, the proof of Theorem 2 is completed by applying \cite[Theorem 6.6]{GT}. The result is based on semi-global practical stability of (\ref{LDsim}) and the recurrence condition on the set $\Delta_e\mathbb{B}$ is needed since the stability property is not global.
\end{proof}
\begin{remark}\label{r7}
	Since the compact set $\Lambda_{\mathrm{opt}}$ is UGAS, it follows from \cite[Theorem 7.21]{GST} that $\Lambda_{\mathrm{opt}}$ is also robustly globally asymptotically stable. Namely, there exists a continuous function $\varsigma:\mathbb{R}^{\bar{n}}\rightarrow\mathbb{R}_{\geq0}$, such that $\Lambda_{\mathrm{opt}}$ is UGAS for $\dot{\tilde{\lambda}}\in-\overline{\text{co}}(\Theta\partial f(\tilde{\lambda}+\varsigma(\tilde{\lambda})\mathbb{B}))+\varsigma(\tilde{\lambda})\mathbb{B}$, where for a set $S$, $\overline{\text{co}}S$ denotes the closed convex hull of $S$. The function $\rho$ can be used to cover some uncertainties commonly seen in applications such as quantization errors, measurement noise. Moreover, a smooth Lyapunov function is guaranteed to exist by a converse theorem given in \cite{GST} which can also be used as a tool for further robustness analysis. \hfill $\blacksquare$
\end{remark}

\begin{remark}\label{round2r2}
  If, in addition to balanced stepsizes \cite[Section 3]{Borkar98}, identical stepsize sequences are used, it can be shown that $\Theta=I$ \cite[Chapter 7]{Borkar1}. In this case one can use $|\lambda|_{\Lambda_{\mathrm{opt}}}^2$ as a Lyapunov function in the proof of Theorem~\ref{mainresult1D} without relying on radial unboundedness of the dual function. 
  \hfill $\blacksquare$
\end{remark}
\section{Additional Results}\label{S-AR}
We first discuss a special case where an almost sure convergence of $\lambda_{[k]}$ in (\ref{reducedmodelfull}) without requiring the iteration to be almost surely bounded is obtained. Then, some convergence rates analysis under more restrictive assumptions are provided. 
%The reason is that one major drawback of the mean differential inclusion analysis is that there is no convenient way for convergence rates analysis of the discrete-time iterations. Thus, we impose more restrictions on the problem to use standard tools from convex optimization to show sublinear convergence rates in a stochastic sense. 
\subsection{Bounded Supergradient}
Almost sure bounded iteration assumption is only used to ensure that $g(\cdot)$ in (\ref{e2.5}) is also almost surely bounded. Thus, it can be established that $\{\bar{\lambda}(t_{[k]}+\cdot)\}$ is equicontinuous and bounded. Thus, if we find another way to ensure the boundedness of $g(\cdot)$, then we can show almost sure convergence without imposing the condition that the iteration is almost surely bounded. One possible approach is to assume that the set $\mathcal{X}$ is compact. In most situations, this is a reasonable assumption since the decision variables are meaningful only when they are finite. The boundedness of the supergradient of $Q(\lambda)$ comes from the fact that the supergradient in the algorithm is formed by the difference between two elements in $\mathcal{X}$ which is bounded.
\begin{corollary}\label{c3}
	Under Assumptions~\ref{a1}--\ref{a6}, if $\mathcal{X}$ is a compact set and $\beta_{[k]}=0$ for all $k\in\mathbb{Z}_{\geq0}$, then the sequence $\{\lambda_{[k]}\}_{k=1}^{\infty}$ generated by (\ref{reducedmodelfull}) satisfies $\lim_{k\rightarrow\infty}|\lambda_{[k]}|_{\Lambda_{\mathrm{opt}}}=0$ almost surely, for any $(\lambda_0,\alpha_0,\gamma_0)\in\mathbb{R}^{\bar{n}}\times\mathbb{R}_{>0}^{|\overline{\mathcal{E}}|}\times\{0\}$.
	\hfill $\blacksquare$
\end{corollary}

\begin{remark}\label{revise4}
  The assumption on the compactness of $\mathcal{X}$ is a simplifying assumption to ensure the almost sure boundedness of the dual variables which is a crucial element of our proof as described in Section~\ref{sec3}. This is in contrast with \cite{TTM} where compactness is to guarantee a differentiable dual function. We do not rely on the differentiability of the dual function in the analysis and the compactness assumption can be replaced by any assumption that results in the boundedness of the dual variable (see \cite{BM} for example). Investigating other sufficient conditions for boundedness of $\lambda$ is left for future work. 
  %In Section~\ref{sec3} we show that convergence of our algorithm requires the dual variable to be almost surely bounded, which can be guaranteed if $\mathcal{X}$ is bounded as shown in Section~\ref{S-AR}. However, unlike \cite{TTM} where compactness is to yield a differentiable dual function, we prove the convergence of our algorithm for case where the dual function may be non-differentiable, which means if alternatives are available to ensure bounded dual updates, compactness of $\mathcal{X}$ will not play a role in our convergence analysis. If $\mathcal{X}$ is bounded, then the diameter of $\mathcal{X}$ can be used as an upper bound of the supergradient of the dual function since it will be shown later in Lemma~\ref{superg} that the supergradient of the dual function used in our updates are obtained via the difference between two primal variables.
	\hfill $\blacksquare$
\end{remark}
\subsection{Convergence Under Stronger Update Assumptions}
In this part, we give another set of sufficient conditions to achieve convergence based on \cite{LSN}. We show that under an assumption on the update rule that is stronger than Assumption~\ref{a4} and global stepsizes\footnote{We require all dual updates to use the same stepsize. It can be local in situations where all updating edges share a common positive end at each time.}, the algorithm almost surely converges if all other assumptions in Theorem~\ref{mainresult1D} are satisfied and there is an appropriate upper bound on the supergradient of the dual function during updates. 
%For easier presentation of the results here, we will use the following reduced model that takes into account the updates for $\lambda$ only since the aim in this part is to show the almost sure boundedness of $\lambda$:
%\begin{equation}\label{e28}
%\lambda^+=\lambda+H(\alpha,x,v^+),
%\end{equation}
%where $v$ here is the placeholder of a random binary variable indicating the edges that are updating. The updates for the stepsize sequences follow the same rule as shown in (\ref{e1.2.3}). As discussed previously, $H_{ij}(\alpha,x,v^+)=\alpha_{(ij)} v_{(ij)}^+(E_{ij}x_{(i)}-E_{ji}x_{(j)})$ can be equivalently written as $H_{ij}(\alpha,x,v^+)\in\alpha_{(ij)} v_{(ij)}^+((\partial Q)_{ij}+e^+)$ where $(\partial Q)_{ij}$ represents the components of the supergradient of $Q$ corresponding to $\lambda_{(ij)}$. 

Instead of Assumption~\ref{a4}, we make the following assumption that can be satisfied if edges are updated in an i.i.d.~fashion:
\begin{assumption}\label{a10}
	Considering system (\ref{reducedmodelfull}), for any $(i,j)\in\overline{\mathcal{E}}$, $\mathbb{P}(v_{(ij),{[k]}}=1|\mathcal{F}_{k-1})=p_{ij}>0$ for all $k\geq1$.
	\hfill $\blacksquare$
\end{assumption}
Note that, Assumption~\ref{a10} is sufficient but not necessary to ensure Assumption~\ref{a4} as it rules out situations where deterministic protocols are used. In addition to Assumption~\ref{a10}, we assume the following upper bound on the supergradient of the dual function $Q$.
\begin{assumption}\label{a3}
	There exists a scalar $c$, such that $
	c^2(1+|\lambda_{[k]}|_{\Lambda_{\mathrm{opt}}}^2)\geq |g(\lambda_{[k]})|^2,$ where $g(\lambda)\in\partial^+ Q(\lambda)$, $\forall k\in\mathbb{Z}_{\geq0},$ 
	almost surely.
	\hfill $\blacksquare$
\end{assumption}

\begin{remark}\label{r8}
% 	The inequality holding for all $\lambda$ is sufficient for the assumption to hold. 
	This assumption is weaker than both the bounded subgradient assumption in \cite{TTM} and the Lipschitz continuous gradient assumption in \cite{NCN}. For example, it can be satisfied if $\mathcal{X}$ is compact. Moreover, if the cost is quadratic which can be commonly seen in linear MPC related applications, then the gradient will be linear which also satisfies Assumption~\ref{a3}.
	\hfill $\blacksquare$
\end{remark}

The following lemma will be used to establish almost sure convergence of $\lambda$ in (\ref{reducedmodelfull}):
\begin{lemma}\cite[Proposition A.4.5]{Bert1}\label{L4}
	Let $\mathcal{F}_k$ be the history of (\ref{reducedmodelfull}) and (\ref{rmcomponent}) up to time $k$. In addition, let $p_{[k]}$, $a_{[k]}$, $w_{[k]}$ and $u_{[k]}$ be sequences of non-negative random variables and let the following relation hold almost surely for any $k\in\mathbb{Z}_{\geq0}$
	\begin{equation}\label{lemm4}
	\mathbb{E}[p_{[k+1]}|\mathcal{F}_k]\leq (1+{a}_{[k]}){p}_{[k]}-{u}_{[k]}+{w}_{[k]},
	\end{equation}
	where $\sum_{k=0}^{\infty}a_{[k]}<\infty$ and $\sum_{k=0}^{\infty}w_{[k]}<\infty$ almost surely. Then almost surely $p_{[k]}$ will converge to some random variable $p$ and $\sum_{k=0}^{\infty}u_{[k]}<\infty$.
	\hfill $\blacksquare$
\end{lemma}
\begin{proposition}\label{prop4}
	Consider Assumptions~\ref{a1},~\ref{a6},~\ref{a10}, and~\ref{a3}. Let $\beta_{[k]}=0$ for all $k\in\mathbb{Z}_{\geq0}$ and all edges use the same stepsize sequence $\{\alpha_{[k]}\}$ satisfying $\sum_{k=0}^{\infty}\alpha_{[k]}=\infty$ and $\sum_{k=0}^{\infty}\alpha_{[k]}^2<\infty$ during dual updates. Then $\lambda_{[k]}$ generated by (\ref{reducedmodelfull}) almost surely converges to a point in $\Lambda_{\mathrm{opt}}$ as $k\rightarrow\infty$ for any $(\lambda_{[0]},\alpha_{[0]},\gamma_{[0]})\in\mathbb{R}^{\bar{n}}\times\mathbb{R}_{>0}^{|\overline{\mathcal{E}}|}\times\{0\}$.
	\hfill $\blacksquare$
\end{proposition}
% \todo[inline]{Any proof that is an exact copy/or a trivial extension of our CDC paper can be removed.}
\begin{proof}
	For any $\lambda^*\in\Lambda_{\mathrm{opt}}$, and $g\in\partial^+ Q$:
	\begin{align*}
	&|\lambda_{[k+1]}-\lambda^*|^2=|\lambda_{[k+1]}-\lambda_{[k]}+\lambda_{[k]}-\lambda^*|^2\\
	&=|\lambda_{[k]}-\lambda^*|^2+|{\alpha}_{[k]}\odot (g(\lambda_{[k]})+e_{[k+1]})\odot v_{[k+1]}|^2\\
	&\ \ \ +2\langle {\alpha}_{[k]}\odot (g(\lambda_{[k]})+e_{[k+1]})\odot v_{[k+1]},\lambda_{[k]}-\lambda^*\rangle.
	\end{align*}
%	\begin{equation*}
%	\begin{split}
%	|\lambda&_{k+1}-\lambda^*|^2=|\lambda_{[k+1]}-\lambda_{[k]}+\lambda_{[k]}-\lambda^*|^2\\
%	&=|\lambda_{[k+1]}-\lambda_{[k]}|^2+|\lambda_{[k]}-\lambda^*|^2\\
%	&\ \ \ +2\langle\lambda_{[k+1]}-\lambda_{[k]},\lambda_{[k]}-\lambda^*\rangle\\
%	&=|\lambda_{[k]}-\lambda^*|^2+|H(\alpha,x,v^+)|^2\\
%	&\ \ \ +2\langle H(\alpha,x,v^+),\lambda_{[k]}-\lambda^*\rangle\\
%	&=|\lambda_{[k]}-\lambda^*|^2+|{\alpha}_{[k]}\odot (s_Q(\lambda_{[k]})+e_{[k+1]})\odot v_{[k+1]}|^2\\
%	&\ \ \ +2\langle {\alpha}_{[k]}\odot (s_Q(\lambda_{[k]})+e_{[k+1]})\odot v_{[k+1]},\lambda_{[k]}-\lambda^*\rangle.
%	\end{split}	  
%	\end{equation*}
  Since $\mathbb{E}[e_{[k+1]}|\mathcal{F}_k]=0$ and
  $\mathbb{E}[|e_{[k+1]}|^2|\mathcal{F}_k]\leq K$, taking expectation of $|\lambda_{[k+1]}-\lambda^*|^2$ conditioned on past history yields
	\begin{align}
	&\mathbb{E}[|\lambda_{[k+1]}-\lambda^*|^2|\mathcal{F}_k]=\mathbb{E}[|{\alpha}_{[k]}\odot e_{[k+1]}\odot v_{[k+1]}|^2|\mathcal{F}_k] \notag \\
	&\ \ \ +|{\lambda}_{[k]}-\lambda^*|^2+\mathbb{E}[|{\alpha}_{[k]}\odot g({\lambda}_{[k]})\odot v_{[k+1]}|^2|\mathcal{F}_k] \notag \\
	&\ \ \ +2\mathbb{E}[\langle {\alpha}_{[k]}\odot g({\lambda}_{[k]})\odot v_{[k+1]},{\lambda}_{[k]}-\lambda^*\rangle|\mathcal{F}_k] \notag \\
	&\leq|{\lambda}_{[k]}-\lambda^*|^2+{\alpha}_{[k]}^2K++{\alpha}_{[k]}^2c^2(1+|{\lambda}_{[k]}-\lambda^*|^2) \notag \\
	&\ \ \ +2{\alpha}_{[k]}\langle Wg({\lambda}_{[k]}),{\lambda}_{[k]}-\lambda^*\rangle \notag \\
	&=(1+{\alpha}_{[k]}^2c^2)|{\lambda}_{[k]}-\lambda^*|^2+{\alpha}_{[k]}^2(K+c^2) \notag \\
	&\ \ \ +2{\alpha}_{[k]}\langle Wg({\lambda}_{[k]}),{\lambda}_{[k]}-\lambda^*\rangle, \label{e.p4.1}
	\end{align}
	where $W:=\text{diag}(\ldots,p_{ij}I,\ldots)\in\mathbb{R}^{\bar{n}}\times\mathbb{R}^{\bar{n}}$ is the block diagonal matrix that contains $p_{ij}$ for all $(i,j)\in\overline{\mathcal{E}}$ and is positive definite due to Assumption~\ref{a10}. Thus, ${W}^{\frac{1}{2}}$ and ${W}^{-\frac{1}{2}}$ are well-defined. Define $m_1:=|W^{-\frac{1}{2}}|^2$, then %almost surely (\ref{e.p4.1}) can be rewritten as:
	\begin{align*}
	&\mathbb{E}[|W^{-\frac{1}{2}}(\lambda_{[k+1]}-\lambda^*)|^2|\mathcal{F}_k]\leq\\
	&(1+{\alpha}_{[k]}^2c^2)|W^{-\frac{1}{2}}({\lambda}_{[k]}-\lambda^*)|^2+{\alpha}_{[k]}^2(|W^{-\frac{1}{2}}|^2|K+c^2)\\
	&+2{\alpha}_{[k]}\langle{W}^{\frac{1}{2}}g(\lambda_{[k]}),{W}^{-\frac{1}{2}}({\lambda}_{[k]}-\lambda^*)\rangle\\
	&\leq(1+{\alpha}_{[k]}^2c^2)|W^{-\frac{1}{2}}({\lambda}_{[k]}-\lambda^*)|^2+{\alpha}_{[k]}^2({m_1K+c^2})\\
	&-2{\alpha}_{[k]}(Q(\lambda^*)-Q(\lambda_{[k]}))  
	\end{align*}
	Since ${\alpha}_{[k]}$ is square summable, we can apply Lemma~\ref{L4} by taking $p_{[k]}=|W^{-\frac{1}{2}}(\lambda_{[k]}-\lambda^*)|^2$, $a_{[k]}={\alpha}_{[k]}^2c^2$, $u_{[k]}=2{\alpha}_{[k]}(Q(\lambda^*)-Q(\lambda_{[k]}))$ and $w_{[k]}={\alpha}_{[k]}^2({m_1K+c^2})$ to arrive at the conclusion that $p_k$ almost surely converges to a random variable which means the algorithm is almost surely bounded. To see this, define $S_A$ as the set of all bounded sequences in a metric space and $S_B$ as the set of all convergent sequences in the same metric space. A convergent sequence in a metric space is bounded. Thus, $S_B\subseteq S_A$. Consider a proper probability measure $P$ such that $P(S_B)\leq P(S_A)$. Since $P(S_B)=1$, we must have $1=P(S_B)\leq P(S_A)\leq 1$ which means $P(S_A)=1$. The rest of the proof follows from the proof of Theorem~\ref{mainresult1D}.
\end{proof}
\subsection{Convergence Rates}
Since we focus on asynchronous algorithms with diminishing stepsizes, one may not expect a convergence rate that is as competitive as those that assume strongly convex functions with Lipschitz gradient. Nevertheless, we are able to state the following result under the hypotheses of Proposition~\ref{prop4}, whose proof follows from the proof of Proposition 1 in \cite{LSN}.
\begin{proposition}\label{Pr1}
	Consider Assumptions~\ref{a1},~\ref{a6},~\ref{a10}, and~\ref{a3}. Let $\beta_{[k]}=0$ for all $k\in\mathbb{Z}_{\geq0}$ and all edges use the same stepsize sequence $\{\alpha_{[k]}\}$ satisfying $\sum_{k=0}^{\infty}\alpha_{[k]}=\infty$ and $\sum_{k=0}^{\infty}\alpha_{[k]}^2<\infty$ during dual updates. Then almost surely $\underset{k\rightarrow\infty}{\liminf}\ k\alpha_{[k]}(Q^*-Q(\lambda_{[k]}))=0$ for any $(\lambda_{[0]},\alpha_{[0]},\gamma_{[0]})\in\mathbb{R}^{\bar{n}}\times\mathbb{R}_{>0}^{|\overline{\mathcal{E}}|}\times\{0\}$. \hfill $\blacksquare$
\end{proposition}
%\begin{proof}
% 	The proof is done by contradiction. Assume that with a positive probability there exists some $\epsilon>0$ and $\bar{k}>0$ such that for all $k\geq \bar{k}$, we have $k\alpha_{[k]}(Q^*-Q(\lambda_{[k]}))\geq\epsilon$.
% 	Then we have $\alpha_{[k]}(Q^*-Q(\lambda_{[k]}))\geq\frac{\epsilon}{k}, \forall k\geq\bar{k}$, which implies
% 	$\sum_{k=\bar{k}}^{\infty}\alpha_{[k]}(Q^*-Q(\lambda_{[k]}))\geq\epsilon\sum_{k=\bar{k}}^{\infty}\frac{1}{k}=\infty$. This contradicts Lemma~\ref{L4} and Proposition~\ref{prop4} and the proof is complete.
%\end{proof}
\begin{remark}\label{r5}
	 After sufficiently many iterations it is possible to observe at least one iteration where the difference between the value obtained at this iteration is arbitrarily close to $Q^*$. %However, it is not guaranteed for iterations that will happen later. 
	 Moreover, the convergence rate depends on how the stepsizes are chosen. Let $\alpha_{[k]}=\frac{1}{(k+1)^q}$, where $\frac{1}{2}<q\leq 1$ for all agents. Define $b_{[k]}:=\underset{i\leq k}{\inf}(Q(\lambda^*)-Q(\lambda_{[k]}))$. If $\underset{k\rightarrow\infty}{\liminf}\ b_k<\epsilon$. Hence, $\underset{k\rightarrow\infty}{\liminf}\ k^{1-q}b_{[k]}=0$. It in turn yields a sublinear convergence rate of $O(\frac{1}{k^{1-q}})$ almost surely. \hfill $\blacksquare$
\end{remark}
%%%%%%%%%%%It is well known \cite{Bert1} that subgradient method with fixed stepsize has a convergence rate of $O(\frac{1}{\sqrt{k}})$ (to the neighborhood of the optimal value). Since $\frac{1}{2}<q\leq1$, the almost sure convergence rate that can be guaranteed is no better than subgradient method. However, this is to be expected due to the fact that we are using diminishing stepsizes to ensure almost convergence to the exact optimal value. It can also be observed that. when $q$ approaches 1, the convergence rate estimates becomes arbitrarily slow. 
We state another sublinear convergence estimate in a stochastic setting under the uniform boundedness of the subgradient at $\lambda_{[k]}$ for all $k\in\mathbb{Z}_{\geq 0}$.
\begin{assumption}\label{a8}
	There exists a scalar $G>0$, such that
	$
	G\geq |g(\lambda_{[k]})|,\ \forall k\in\mathbb{Z}_{\geq0}, g(\lambda)\in\partial^+ Q(\lambda)
	$
	almost surely.
	\hfill $\blacksquare$
\end{assumption}

\begin{proposition}\label{T3}
	Consider the sequence of random variables $\lambda_{[k]}$ generated by (\ref{reducedmodelfull}). Under Assumptions~\ref{a1},~\ref{a6},~\ref{a10}, and~\ref{a3}, assume $\beta_{[k]}=0$ for all $k\in\mathbb{Z}_{\geq0}$ and all edges use the same stepsize sequence $\{\alpha_{[k]}\}$ satisfying $\sum_{k=0}^{\infty}\alpha_{[k]}=\infty$ and $\sum_{k=0}^{\infty}\alpha_{[k]}^2<\infty$ during dual updates. Then for any $(\lambda_{[0]},\alpha_{[0]},\gamma_{[0]})\in\mathbb{R}^{\bar{n}}\times\mathbb{R}_{>0}^{|\overline{\mathcal{E}}|}\times\{0\}$: 
	\begin{equation*}
	Q^*-Q_{\text{best}[k]}\leq \frac{m_1R^2+(m_1K+G^2)\sum_{i=1}^{k}\alpha^2_{[i]}}{2\sum_{i=0}^{k}\alpha_{[i]}},
	\end{equation*}
	almost surely, where $Q_{{\text{best}}[k]}=\underset{i=0,1,...,k}{\text{max}}\mathbb{E}[Q(\lambda_{[k]})|\mathcal{F}_{i-1}]$ and $R:=\underset{\lambda^*\in \Lambda_{\mathrm{opt}}}{\text{sup}}|\lambda_{[0]}-\lambda^*|$. \hfill $\blacksquare$
\end{proposition}
\begin{proof}
  Similar to the proof of Proposition~\ref{prop4}:
	\begin{align}
	\mathbb{E}[&|W^{-\frac{1}{2}}(\lambda_{[k+1]}-\lambda^*)|^2|\mathcal{F}_k]\leq|W^{-\frac{1}{2}}({\lambda}_{[k]}-\lambda^*)|^2\notag\\
	&+\alpha_{[k]}^2(m_1K+G^2)-2{\alpha}_{[k]}(Q(\lambda^*)-Q(\lambda_{[k]})) \label{e.p6.1}
 	\end{align}
 	The rest of the proof follows from the proof of Theorem 2 in \cite{LSN} and is omitted for the sake of brevity.
% 	Similar to the proof of Proposition~\ref{prop4}:
% 	\begin{align}
% 	\mathbb{E}[&|W^{-\frac{1}{2}}(\lambda_{[k+1]}-\lambda^*)|^2|\mathcal{F}_k]\leq|W^{-\frac{1}{2}}({\lambda}_{[k]}-\lambda^*)|^2\notag\\
% 	&+\alpha_{[k]}^2(m_1K+G^2)-2{\alpha}_{[k]}(Q(\lambda^*)-Q(\lambda_{[k]})) \label{e.p6.1}
% 	\end{align}
% 	Applying the law of total expectations iteratively to (\ref{e.p6.1}) \cite[Proposition D.5 (b), Appendix D]{BT}, we have
% 	\begin{align}
% 	\mathbb{E}[|W^{-\frac{1}{2}}&(\lambda_{[k+1]}-\lambda^*)|^2|\mathcal{F}_k]\leq \mathbb{E}[|W^{-\frac{1}{2}}({\lambda}_{[k]}-\lambda^*)|^2|\mathcal{F}_{k-1}]\notag\\
% 	&+\alpha_{[k]}^2(m_1K+G^2)-2\alpha_{[k]}(Q^*-\mathbb{E}[Q(\lambda_{[k]})|\mathcal{F}_{k-1}])\notag\\
% 	&\leq |W^{-\frac{1}{2}}(\lambda-\lambda_{[0]})|^2+\sum_{i=0}^{k}\alpha^2_{[i]}(m_1K+G^2) \notag\\
% 	&-2\sum_{i=0}^{k}\alpha_{[i]}(Q^*-\mathbb{E}[Q(\lambda_{[i]})|\mathcal{F}_{i-1}]), \label{eT2}
% 	\end{align}
% 	almost surely. Since $\underset{\lambda^*\in \Lambda_{\mathrm{opt}}}{\text{sup}}|\lambda_{[0]}-\lambda^*|=R$, it follows that $0\leq m_1{R^2}-2\sum_{i=0}^{k}\alpha_{[i]}(Q^*-\mathbb{E}[Q(\lambda_{[i]})|\mathcal{F}_{i-1}])+\sum_{i=0}^{k}\alpha^2_{[i]}(m_1K+G^2)$. Then, almost surely we have
% 	\begin{equation*}
% 	Q^*-Q_{\text{best}[k]}\leq \frac{m_1R^2+(m_1K+G^2)\sum_{i=1}^{k}\alpha^2_{[i]}}{2\sum_{i=0}^{k}\alpha_{[i]}}.
% 	\end{equation*}
% 	Since $\sum_{i=0}^{k}\alpha_{[i]}=\infty\ \text{and}\ \sum_{i=0}^{k}\alpha_{[i]}^2<\infty$, we establish almost sure convergence of $Q_{\text{best}[k]}$ to $Q^*$.
\end{proof}
% \begin{remark}
% 	Based on Proposition~\ref{T3}, we can state an improved convergence rate. For the widely used stepsizes $\alpha_{[k]}=\frac{1}{k+1}$ compared to the one established in Proposition~\ref{Pr1}. By using the fact that $\sum_{t=0}^{k}\frac{1}{t+1}\geq\int_{0}^{k}\frac{1}{1+t}dt=\log (t+1)|_{0}^{k}=\log (k+1)$, a convergence rate (in a different setting) of $O(\frac{1}{\log k})$ for this choice of stepsize is established. 
% % 	Convergence rates for stepsizes of the form $\alpha_{[k]}=\frac{1}{(k+1)^q}$, $\frac{1}{2}<q\leq 1$ can also be stated using the same approach but they happen to have the same order as the one shown in Proposition~\ref{Pr1}.
% 	\hfill $\blacksquare$
% \end{remark}
\begin{remark}
% 	Based on Proposition~\ref{T3}, we can state an improved convergence rate. For the widely used stepsizes $\alpha_{[k]}=\frac{1}{k+1}$ compared to the one established in Proposition~\ref{Pr1}. 
	Since $\sum_{t=0}^{k}\frac{1}{t+1}\geq\int_{0}^{k}\frac{1}{1+t}dt=\log (t+1)|_{0}^{k}=\log (k+1)$, a convergence rate of $O(\frac{1}{\log k})$ for $\alpha_{[k]}=\frac{1}{k+1}$ is established using Proposition~\ref{T3}. 
% 	Convergence rates for stepsizes of the form $\alpha_{[k]}=\frac{1}{(k+1)^q}$, $\frac{1}{2}<q\leq 1$ can also be stated using the same approach but they happen to have the same order as the one shown in Proposition~\ref{Pr1}.
	\hfill $\blacksquare$
\end{remark}

\begin{remark}\label{revise6}
	Proposition~\ref{T3} establishes the sublinear convergence of the dual algorithm. If $F_i$ in \eqref{PP} is strongly convex and has Lipschitz gradient, then it can be shown via Fenchel duality that the dual function $Q$ is strongly concave and has Lipschitz continuous gradient. Then, if we further assume that i.i.d.~local timers are used in the updates, we recover the random block coordinate descent algorithm of \cite[Theorem 1]{Wright} where a stochastic linear convergence rate can be achieved. As a result, while the main focus of this work is the asynchronous aspect of the algorithm rather than achieving the fastest possible convergence rate, linear convergence can be recovered if we make more stringent assumptions about the problem and update rules in the network.   
	\hfill $\blacksquare$
\end{remark}

\section{A Numerical Example}\label{S-NE}
We consider the problem of minimizing $F(x):=\sum_{i=1}^{50}F_i(x)$ over a network of 50 agents. In this case, all selection matrices $E_{ij}$ are identity matrices and for simplicity we choose $e$ to be $0$ so that the superdradient of the dual is error free. The functions $F_i(x)$ are given by $F_i(x)=\max\{-w_{i}(x-a_i)+b_s,b_s\}$ (from equation (5) in \cite{HI}) for $1\leq i\leq 5$, where $b_s$ is selected to be $0$ while $w_{i}$ and $a_i>0$ are randomly chosen in the intervals $[0.2, 1]$ and $[2,8]$ respectively. Note that these functions are non-smooth convex but not strictly convex. For $6\leq i\leq 50$, $F_i(x)=x\log p_ix$, where for each $1\leq i\leq50$, $p_i$ is a real number chosen in $(1,5)$ uniformly at random. We consider two communication graphs: a path graph and a connected random geometric graph with 358 edges with all dual variables initialized at $0$.
\begin{figure}%
	\centering
	\subfloat[Plot of $|Q(\lambda^*)-Q(\lambda)|$ versus the number of iterations]{{\includegraphics[width=.45\textwidth]{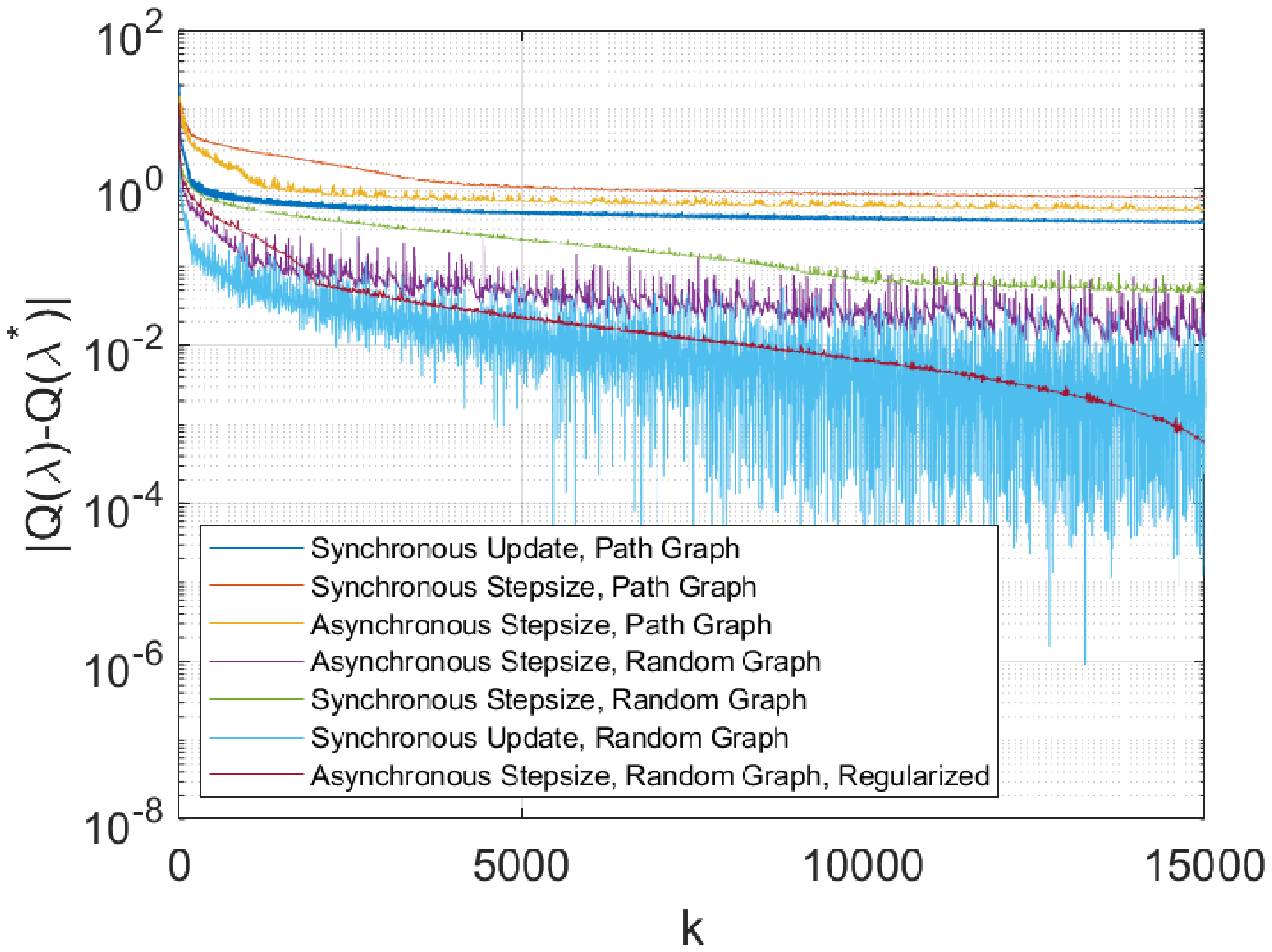} }}%
	\ 
	\subfloat[Plot of primal variables $x_{(i)}$ versus the number of iterations]{{\includegraphics[width=.45\textwidth]{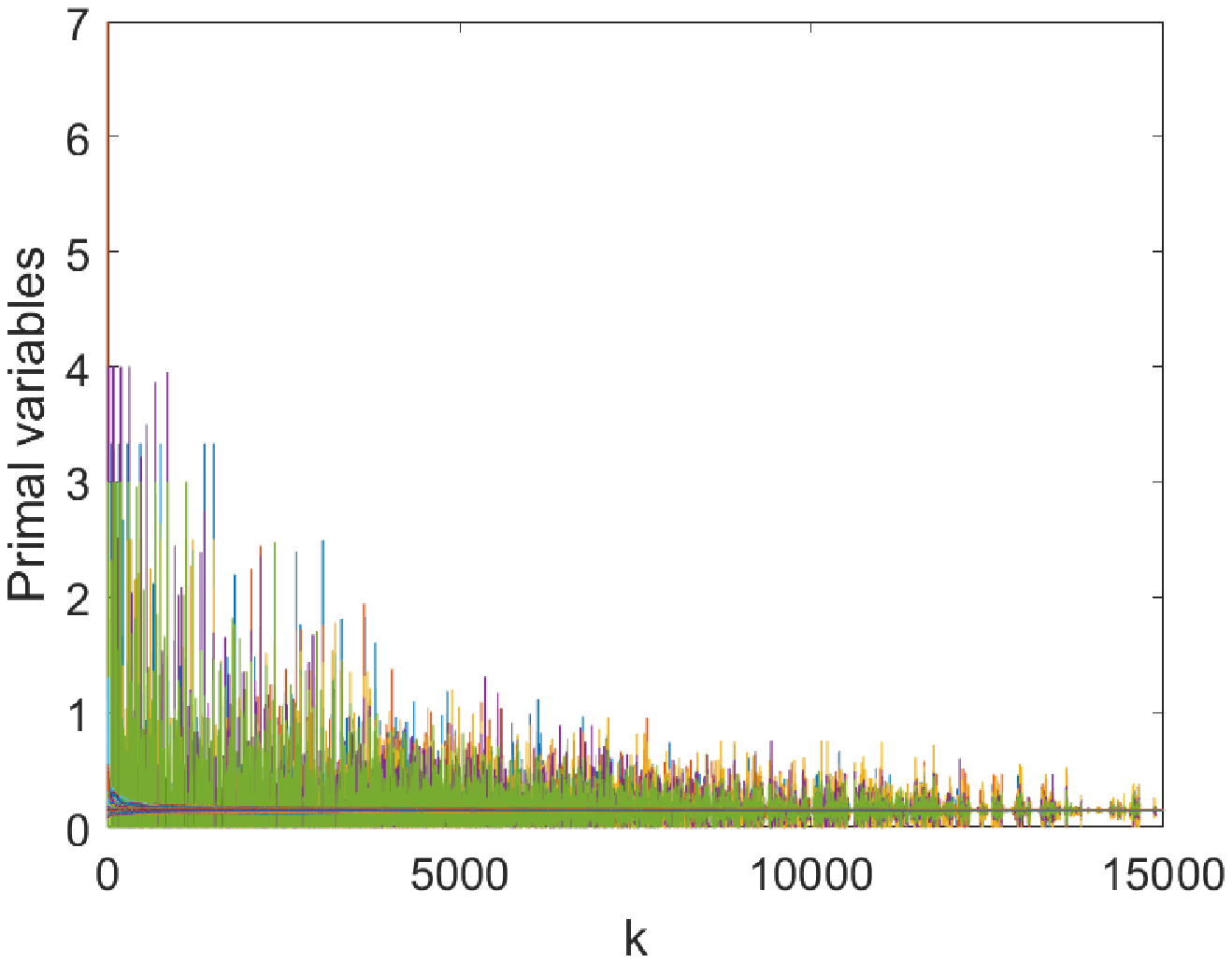} }}%
	\qquad
	\caption{Fig. 1(a) shows the convergence of the dual function under different settings, while Fig. 1(b) shows the primal variables indeed approach consensus as $k$ grows. }%
	\label{figdual}%
\end{figure}
%\begin{figure}[h]
%	\centering
%	\subfloat[Path graph with 50 nodes and 98 edges]{{\includegraphics[width=3.9cm,height=3cm]{ring_network.eps} }}%
%	\qquad
%	\subfloat[Randomly generated graph with 50 nodes and 358 edges]{{\includegraphics[width=3.9cm,height=3cm]{random_network.eps} }}%
%	\caption{Two communication graphs}%
%	\label{commgraph}%
%\end{figure}
Each agent $i\in\mathcal{V}$ is given a local decision variable $x_{(i)}$. Since for $1\leq i\leq5$, $F_i(x)$ is not strictly convex, as discussed in Remark~\ref{revise2}, only convergence of the dual variable is guaranteed. We consider the following update scheme. We assign to each $i\in\mathcal{V}$ a real number $\tilde{p}_i$ chosen in $(0.5,1)$ uniformly at random, a local counter $\varpi_i$ that counts the number of times that there is at least one update at edges connecting $i$ to its neighbors in the last 10 iterations in the whole network. Then for $(i,j)\in\overline{\mathcal{E}}$, $\lambda_{(ij)}$ gets updated with a probability $\tilde{p}_i\tilde{p}_j0.7^{\varpi_i\varpi_j}$. The probability of updating for each node $i\in\mathcal{I}$ reduces as the number of updates grows in the last 10 global iterations. It is worth mentioning that although $\varpi_i$ uses global information in the network, the introduction of $\varpi_i$ is merely to simulate the algorithm with time varying updating probabilities. The implementation of the algorithm does not require global information in the network as long as the assumptions of the main results are satisfied. The stepsizes for dual updates are assumed to be the same, initialized at 0.15, and generated by $A(\alpha) = 0.15(1+ \log_{0.51}0.15 - \log_{0.51}\alpha )^{-0.51}$ (This stepsize rule corresponds to $\alpha_{[k]}= 0.15(1+k)^{-0.51}$ in the case of a global time counter $k$). For each of the two networks, we simulate the algorithm for the case where all communication and updates are synchronous, the updates are asynchronous with global stepsize, and the updates are asynchronous with local stepsize. To compare, we also consider the case where for $1\leq i\leq5$, $F_i(x)$ is regularized by $0.005x^2$ so that we should have all $x_{(i)}$, $i\in\mathcal{V}$ converge to the same $x^*$ that minimizes $F(x)$. Dual variables are all initialized at the zero vector. Furthermore, since we work with non-strictly convex functions, to satisfy Assumption~\ref{a1}, we work with $x$ over compact sets. By doing so, we can ensure a real solution to the local problem always exists and the upper bound of $x$ does not need to be known, The plots for dual functions and primal variables when regularizers are used are given in Fig.~\ref{figdual}.

As predicted by our results, $Q(\lambda)$ tends to its maximum and all agents reach consensus as the number of iterations tends to infinity if regularizers are used. Moreover, the use of local stepsizes are observed to provide better convergence rates. This can be explained by the fact that the stepsize for the edges that are not updating does not decrease while in the global stepsize case the stepsize decreases whenever there is an update anywhere in the network. Although not theoretically explored in the paper, adding more communication links in the network is also observed to accelerate the convergence of both dual and primal variables which may also induce more fluctuations caused by the numerical errors. It can be seen from Fig.~\ref{figdual}(a) that the three dual functions corresponding to the random graph converge faster than the three dual functions corresponding to the path graph but with more fluctuations. It shows an interesting trade-off between the number of communication links and a faster convergence rate. The investigation of this interesting relationship between communication and convergence is left as future work. On the other hand, Fig.~\ref{figdual}(b) shows that the primal variables reach consensus asymptotically as the functions are regularized.  
%%%%%%%%%%%%%%%%%%%%%%%%%%%%%%%%%%%%%%%%%%%%%%%%%%%%%%%%%%%%%%%%%%%%%%%%%%%%%%%%%%%%%%%%%%%%%%%%%%%%%%%
\section{Conclusion and Future Work}\label{S-FW}
We have proposed an edge based asynchronous distributed dual algorithm to minimize the sum of convex cost functions that have partially overlapping dependences. The analysis was done by viewing the asynchronous algorithm as block coordinate supergradient update for a concave function and applying techniques from stochastic approximation literature. Under relatively weak assumptions on the cost function, sufficient conditions on the communications between agents and stepsizes were provided to prove almost sure convergence of the algorithm to the solution of the dual problem. We gave a sublinear convergence rate estimate for the algorithm under stronger assumptions. A numerical example was given to illustrate our main result. Future work will focus on providing convergence rate analysis under weaker assumptions on the communication between agents and the use of constant stepsizes in the algorithm to achieve faster convergence.

\appendix
\subsection{Background Definitions and Lemmas}
For any $T>0$ and $n\in\mathbb{Z}_{>0}$, the space of continuous functions from $[0,T]$ to $\mathbb{R}^n$ is denoted by $C([0,T];\mathbb{R}^n)$ and the space of measurable functions $f$ from $[0,T]$ to $\mathbb{R}^n$ satisfying $\int_{0}^{T}|f(t)|^2dt<\infty$ is denoted by $\mathcal{L}_2([0,T];\mathbb{R}^n)$. A subset of a topological space is \textit{relatively compact} if its closure is compact. A set $B\subseteq C([0,T];\mathbb{R}^n)$ is \textit{equicontinuous at} $t\in[0,T]$ if for any $\epsilon>0$, there exists a $\delta>0$ such that $|t-s|<\delta$, $s\in[0,T]$ implies $\sup_{f\in B}||f(t)-f(s)||_{\infty,T}<\epsilon$, where $||f||_{\infty,T}:=\sup_{s\in[0,T]}|f(s)|$. It is \textit{equicontinuous} if it is equicontinuous at all $t\in[0,T]$. It is \textit{pointwise bounded} if for any $t\in[0,T]$, $\sup_{f\in B}||f(t)||_{\infty,T}<\infty$. The following result known as the Arzel\`{a}-Ascoli Theorem characterises relative compactness in $C([0,T];\mathbb{R}^n)$.
\begin{lemma}\cite[Appendix A, Theorem 1]{Borkar1}\label{AA}
	A set $B\subseteq C([0,T];\mathbb{R}^n)$ is relatively compact if and only if it is equicontinuous and pointwise bounded. \hfill $\blacksquare$
\end{lemma}
Consider the Hilbert space $\mathcal{L}_2([0,T];\mathbb{R}^n)$ with the inner product defined by $\langle f, g\rangle_T:=\int_{0}^{T}\langle f(t),g(t)\rangle dt,\ f,g\in\mathcal{L}_2([0,T];\mathbb{R}^n)$ and the norm $||f||_{2,T}:=\sqrt{\langle f,f\rangle_T}$. The \textit{weak} topology on $\mathcal{L}_2([0,T];\mathbb{R}^n)$ is defined as the coarsest topology
with respect to which the functions $f\rightarrow\langle f,g\rangle_T$ are continuous for all $g\in\mathcal{L}_2([0,T];\mathbb{R}^n)$. Then, we say $f_n\rightarrow f$ weakly in $\mathcal{L}_2([0,T];\mathbb{R}^n)$ if and only if $\langle f_n,g\rangle_T\rightarrow\langle f,g\rangle_T$ for all $g\in\mathcal{L}_2([0,T];\mathbb{R}^n)$. The following result which is a special instance of the Banach-Alaoglu Theorem characterises relative compactness in $\mathcal{L}_2([0,T];\mathbb{R}^n)$.
\begin{lemma}\cite[Appendix A, Theorem 3]{Borkar1}\label{BA}
	A $||\cdot||_{2,T}$-bounded set $B\subseteq \mathcal{L}_2([0,T];\mathbb{R}^n)$ is relatively compact in the weak topology. \hfill $\blacksquare$
\end{lemma}
Next, we state a martingale convergence theorem that is critical in the proof of Theorem~\ref{mainresult1D}. Let $(\Omega,\mathcal{F},\mathbb{P})$ be a probability space and $\{\mathcal{F}_n\}$ be a family of increasing $\sigma$-subfields of $\mathcal{F}$. For a martingale $\{X_n\}$ with respect to $\{\mathcal{F}_n\}$, we have the following convergence theorem.
\begin{lemma}\cite[Appendix C, Theorem 11]{Borkar1}\label{MC}
	If $\mathbb{E}[|X_n|^2]<\infty$ for any $n\in\mathbb{Z}_{\geq0}$ and $\sum_{n}\mathbb{E}[|X_{n+1}-X_n|^2|\mathcal{F}_n]<\infty$, then $\{X_n\}$ converges almost surely. \hfill $\blacksquare$
\end{lemma}
Some useful properties of convex functions are stated below.
\begin{lemma}\label{LL}
	Let $f:\mathbb{R}^n\rightarrow\mathbb{R}$ be a real-valued convex function. Then, the following statements hold for any $x\in\mathbb{R}^n$:
	\begin{enumerate}[(i)]
		\item $\partial f(x)$ is non-empty, compact and convex.
		\item $f$ is locally Lipschitz.
		\item $\partial f(x)$ is locally bounded.
		\item $\partial f(x)$ is upper semicontinuous\footnote{A set-valued mapping $M:\mathbb{R}^n\rightrightarrows\mathbb{R}^n$ is upper semicontinuous if for each $x\in\mathbb{R}^n$ and $\varepsilon>0$ there exists $\delta>0$ such that $M(x+\delta\mathbb{B})\subseteq M(x)+\varepsilon\mathbb{B}$.}.
		\item $\partial f(x)$ is outer semicontinuous.\hfill $\blacksquare$
	\end{enumerate} 
\end{lemma}
\begin{proof}
	Item (i) follows from \cite[Proposition 3.1.1]{Bert1}. Items (ii) and (iii) follow from \cite[Proposition 3.1.2]{Bert1}. By \cite[Proposition 2.2.7]{Clarke}, if $f(x)$ is convex, then $\partial f(x)$ coincides with the generalized gradient which is upper semicontinuous by \cite[item(d), Proposition 2.1.5]{Clarke}. Lastly, by \cite[Lemma 5.15]{GST}, when $\partial f(x)$ is closed and locally bounded, upper semicontinuity of $\partial f(x)$ is equivalent to outer semicontinuity of $\partial f(x)$. Since $\partial f(x)$ is indeed closed and locally bounded by item (i) and (iii) of Lemma~\ref{LL}, item (v) follows.
\end{proof}
 A function $\psi:\mathbb{R}_{\geq0}\rightarrow\mathbb{R}_{\geq0}$ is \textit{of class} $\mathcal{K}$ if it is
	continuous, strictly increasing and $\psi(0)=0$. It is \textit{of class} $\mathcal{K}_\infty$ if it is
	of class $\mathcal{K}$ and unbounded. A function $\chi:\mathbb{R}_{\geq0}\rightarrow\mathbb{R}_{\geq0}$ is \textit{of class} $\mathcal{L}$ if it is non-increasing and $\lim_{t\rightarrow\infty}\chi(t)=0$. A function $\beta:\mathbb{R}_{\geq0}\times\mathbb{R}_{\geq0}\rightarrow\mathbb{R}_{\geq0}$ is \textit{of class} $\mathcal{KL}$ if $\beta(\cdot,t)\in\mathcal{K}$ for each $t\in\mathbb{R}_{\geq0}$ and $\beta(s,\cdot)\in\mathcal{L}$ for each $s\in\mathbb{R}_{\geq0}$
	Let $V^\circ(x;v):=\limsup_{h\rightarrow 0^+,y\rightarrow x}\frac{V(y+hv)-V(y)}{h}$ be the generalized directional derivative of Clarke of a locally Lipschitz function $V$ at $x$ in the direction of $v$. This is the standard directional derivative $\langle \nabla V(x),v \rangle$ for continuously differentiable $V$ \cite{Clarke}.
	
\subsection{Proof of Lemma~\ref{newstepsize}}
	Since $\bar{\alpha}=\underset{i\in \mathcal{I}_a}{\max}\ \alpha_{(i)}$, all elements of the sequence $\{\bar{{\alpha}}_{[k]}\}$ must belong to the set $\mathcal{M}:=\{a\in\mathbb{R}:a=\alpha_{(i),[k]},\ i\in\mathbb{Z}_{|\overline{\mathcal{E}}|}\backslash\{0\},k\in\mathbb{Z}_{\geq0}\}$, which is the set of all elements of local stepsize sequences. Thus, the sum of squares of $\{\bar{{\alpha}}_{[k]}\}$ must be upper bounded by the sum of squares of all elements in $\mathcal{M}$. Note that, for any $i$ in the finite set $\mathbb{Z}_{|\overline{\mathcal{E}}|}\backslash\{0\}$, the sequence $\{\alpha_{(i),[k]}\}$ is square summable. Thus, all elements in $\mathcal{M}$ must also be square summable. That is $\sum_{k=0}^{\infty}\bar{{\alpha}}_{[k]}^2<\sum_{i=1}^{|\overline{\mathcal{E}}|}\sum_{k=0}^{\infty}\alpha_{(i),[k]}^2<\infty$ for $\{\alpha_{(i),[k]}\}$ generated by $\alpha_{(i)}^+=A_{(i)}(\alpha_{(i)})$. By Assumption~\ref{a4}, $\mathcal{I}_a\neq\emptyset$ infinitely often, almost surely. Therefore, $\bar{{\alpha}}_{[k]}>0$ for infinitely many $k$, almost surely. Since $|\overline{\mathcal{E}}|$ is finite, $\bar{{\alpha}}_{[k]}>0$ for infinitely many $k$ implies that $\{\bar{{\alpha}}_{[k]}\}$ must contain infinitely many terms from at least one sequence $\{{\alpha}_{(i),[k]}\}$ for some $i\in\mathbb{Z}_{|\overline{\mathcal{E}}|}\backslash\{0\}$. Otherwise $\{{\alpha}_{(i),[k]}\}$ will only have a finite number of elements which contradicts the fact that $\bar{{\alpha}}_{[k]}>0$ for infinitely many $k$. Therefore $\sum_{k=0}^{\infty}\bar{{\alpha}}_{[k]}=\infty$ almost surely.
\subsection{Proof of Proposition~\ref{martingaleconvergence}}
	By Lemma~\ref{newstepsize}, $\{\bar{{\alpha}}_{[k]}\}_{k=1}^{\infty}$ is square summable almost surely and the largest element of $\theta(t)$ is $1$ for all $t$. Hence, $\sum_{k=1}^{\infty}\mathbb{E}[|{\zeta}_{[k+1]}-{\zeta}_{[k]}|^2 |\mathcal{F}_k]\leq\sum_{k=1}^{\infty}\bar{{\alpha}}^2_{[k]}\mathbb{E}[|e_{[k+1]}|^2|\mathcal{F}_k]\leq\sum_{k=1}^{\infty}\bar{{\alpha}}^2_{[k]}K<\infty$ almost surely, due to Assumption~\ref{a6} and Lemma~\ref{newstepsize}. Since $\beta_{[k]}=0$ for all $k$, $\mathbb{E}[\langle\bar{{\alpha}}_{[m_1]}\theta(t_{[m_1]})e_{[{m_1}+1]},\bar{{\alpha}}_{[m_2]}\theta(t_{[m_2]})e_{[{m_2}+1]}\rangle]=0$ holds for any $m_1\neq m_2\in\mathbb{Z}_{\geq0}$. Thus, for any $k\geq1$, $\mathbb{E}[|\zeta_{[k]}|^2]=\sum_{k=1}^{\infty}\mathbb{E}[|\bar{{\alpha}}_{[k-1]}e_{[k]}|^2]<\infty$ follows from the fact that $|\theta(t)|$ is bounded. Consequently, all assumptions of Lemma~\ref{MC} are satisfied and ${\zeta}_{[k]}$ converges almost surely to a random variable. The definition of $r(t)$ and Lemma~\ref{newstepsize} yield $r(t_{[k]})\rightarrow\infty$ as $k\rightarrow\infty$. Therefore, the sequence $\{{\zeta}_{[r(t_{[k]}+\cdot)]}-{\zeta}_{[r(t_{[k]})]}\}_{k=1}^{\infty}$ converges to $0$ uniformly on $[0,T]$ for any $T>0$ almost surely when $k$ is sufficiently large. Thus the proof is complete.
\subsection{Proof of Proposition~\ref{lemma3}}
	From Proposition~\ref{martingaleconvergence}, the sequence $\{{\zeta}_{[r(t_{[k]}+u)]}-{\zeta}_{[r(t_{[k]})]}\}_{k=1}^{\infty}$ almost surely converges to $0$ on $[0,\infty)$. Moreover, due to Lemma~\ref{newstepsize}, $t_{[k+1]}-t_{[k]}=\bar{{\alpha}}_{[k]}\rightarrow0$ as $k\rightarrow\infty$, thus for $t\in[t_{[k]},t_{[k+1]})$, $\int_{0}^{u}\theta(t+s)g(\lambda_{[r(t+s)]})ds-\int_{0}^{u}\theta(t+s)g(\bar{\lambda}(t+s))ds\rightarrow0$ as $k\rightarrow\infty$ and the result follows since $k\rightarrow\infty$ implies $t_{[k]}\rightarrow\infty$ by Lemma~\ref{newstepsize}.
\subsection{Proof of Proposition~\ref{Prop2}}
	The proof is based on \cite[Chapter 7, Theorem 2]{Borkar1}. First, let $\tilde{\theta}(t):=(\tau_{(1),[k]},\ldots,\tau_{(|\overline{\mathcal{E}}|),[k]}), \text{for}\ t\in[t_{[k]},t_{[k+1]})$ which is the vector containing all diagonal elements of $\theta(t)$. Then, we view $\tilde{\theta}(\cdot)$ as an element of $\mathcal{V}$, where $\mathcal{V}$ is the space of measurable maps $y(\cdot):[0,\infty)\rightarrow[0,1]^{\bar{n}}$ with the coarsest topology that renders the maps
	$y(\cdot)\rightarrow\int_{0}^{T}\langle g(s),y(s)\rangle ds$ continuous for all $T>0,g(\cdot)\in\mathcal{L}_2([0,T];\mathbb{R}^{\bar{n}})$. This defines the weak topology on $\mathcal{L}_2([0,T];\mathbb{R}^{\bar{n}})$ (see Lemma~\ref{BA}). Since $g(\cdot)\in\mathcal{L}_2([0,T];\mathbb{R}^{\bar{n}})$ and all elements of $\tilde{\theta}(t)$ are uniformly bounded for all $t\geq0$, $\mathcal V$ is a bounded subset of $\mathcal{L}_2([0,T];\mathbb{R}^{\bar{n}})$. Then by the Banach-Alaoglu Theorem (see Lemma~\ref{BA}), $\mathcal{V}$ is relatively compact. Thus, by choosing an appropriate common subsequence of $\{\tilde{\theta}(t_{[k]}+\cdot)\}$ and $\{\bar{\lambda}(t_{[k]}+\cdot)\}$ indexed by $\{k'\}$, we can show that any limit point of $\{\bar{\lambda}(t_{[k]}+\cdot)\}$ in $C([0,T];\mathbb{R}^{\bar{n}})$ denoted by $\tilde{\lambda}(\cdot)$ is of the form
	\begin{equation*}
	\lim_{k'\rightarrow\infty}\bar{\lambda}(t_{[k']}+t):=\tilde{\lambda}(t)=\tilde{\lambda}(0)-\int_{0}^{t}\Theta(s)g(\tilde{\lambda}(s))ds,
	\end{equation*} 
	where $\Theta(\cdot)=\text{diag}(\tilde{\theta}^*_{(1)}(\cdot),\ldots,\tilde{\theta}^*_{(|\overline{\mathcal{E}}|)}(\cdot))$ and $\tilde{\theta}^*(\cdot)=(\tilde{\theta}^*_{(1)}(\cdot),\ldots,\tilde{\theta}^*_{(|\overline{\mathcal{E}}|)}(\cdot))$ is the subsequential limits of $\{\tilde{\theta}(t_{[k']}+\cdot)\}$ in $\mathcal{V}$ corresponding to $\{k'\}$. The limit point of $\{g(\bar{\lambda}(t_{[k']}+\cdot))\}$ in $\mathcal{L}_2([0,T];\mathbb{R}^{\bar{n}})$ belongs to $\partial f(\tilde{\lambda}(\cdot))$ since $\partial f$ is outer semicontinuous with closed and convex values which is guaranteed by Lemma~\ref{LL}.
	%	\footnote{Detailed derivations can be found in the proof of Theorem 4.2 in \cite{BHS}.}. 
	
	For any given $n\geq0$ and $s>0$, let $N(n,s)$ be the smallest integer $m$ such that $t_{[m]}\geq t_{[n]}+s$. Since $\tilde{\theta}(\cdot)$ is a piecewise constant function by definition and $\bar{{\alpha}}_{[k]}=t_{[k+1]}-t_{[k]}$, for any given $i$ and $k$, $\int_{t_{[k]}}^{t_{[k+1]}}\tilde{\theta}_{(i)}(t)dt=\frac{\alpha_{(i),[\gamma_{(i),[k]}]}v_{(i),[k+1]}}{\bar{\alpha}_{[k]}}\bar{\alpha}_{[k]}$. As a result, the integral of $\tilde{\theta}(\cdot)$ can be replaced by appropriate sums. Then for any $i\in\mathbb{Z}_{|\overline{\mathcal{E}}|}\backslash\{0\}$, $\int_{t}^{t+s}\tilde{\theta}^*_{(i)}(y)dy=\lim_{k'\rightarrow\infty}\sum_{m=k'}^{N(k',s)}\frac{\alpha_{(i),[\gamma_{(i),[m]}]}v_{(i),[m+1]}}{\bar{\alpha}_{[m]}}\bar{\alpha}_{[m]}$ which is equal to $\lim_{k'\rightarrow\infty}\sum_{m=\gamma_{(i),[k']}}^{\gamma_{(i),[N(k',s)]}}\alpha_{(i),[m]}$. The last equality follows from the fact that $\alpha_{(i)}$ and $\gamma_{(i)}$ get updated if and only if $v_{(i)}^+=1$ which means at time $k$, the stepsize $\alpha_{(1),[k]}$ generated by (\ref{rmcomponent}) is the same as $\alpha_{(i),[\gamma_{(i),[k]}]}$ generated from the pre-selected update rule $\alpha_{(i)}^+=A_{(i)}(\alpha_{(i)})$. By Assumptions~\ref{a4},~\ref{a5}, $\frac{\int_{t}^{t+s}\tilde{\theta}^*_{(i)}(y)dy}{\int_{t}^{t+s}\tilde{\theta}^*_{(j)}(y)dy}=\lim_{k'\rightarrow\infty}\frac{\sum_{m=\gamma_{(i),[k']}}^{\gamma_{(i),[N(k',s)]}}\alpha_{(i),[m]}}{\sum_{m=\gamma_{(j),[k']}}^{\gamma_{(j),[N(k',s)]}}\alpha_{(j),[m]}}\geq\liminf_{k'\rightarrow\infty}\frac{\Delta\alpha_{(i),[\gamma_{(i),[N(k',s)]}]}}{\alpha_{(j),[\gamma_{(j),[k']}]}}\geq\min\{\Delta\bar{\delta},\frac{\Delta}{\bar{\Delta}}\}:=\varepsilon$ for any $i,j\in\mathbb{Z}_{|\overline{\mathcal{E}}|}\backslash\{0\}$. The first inequality follows from replacing the sum in the numerator (respectively denominator) by the sum of its smallest (resp. largest) element (assuming without loss of generality all stepsize sequences are non-increasing) and Assumption~\ref{a4} that the ratio of the number of updates is lower bounded by $\Delta$. Since $s$ is arbitrary and all diagonal elements of $\Theta(t)$ are upper bounded by $1$ by construction, $\Theta(t)$ has diagonal entries in $[\varepsilon,1]$.
\subsection{Proof of Proposition~\ref{limi}}
	By Proposition~\ref{Prop2}, $\Theta(t)\in\mathbb{R}^{\bar{n}\times \bar{n}}$ is a diagonal matrix with diagonal elements in the interval $[\varepsilon,1]$ for all $t$. Since $f(\tilde{\lambda})$ is locally Lipschitz, the radially unbounded Lyapunov function $V(\tilde{\lambda}):=f(\tilde{\lambda})+Q^*$ is also locally Lipschitz. Moreover, by convexity of $f$, both $V$ and $\nabla f$ vanish on $\Lambda_{\mathrm{opt}}$ and nowhere else. Define $\phi_1(s):=\inf_{|\tilde{\lambda}|_{\Lambda_{\mathrm{opt}}}\geq s}V(\tilde{\lambda})$ and $\phi_2(s):=\sup_{|\tilde{\lambda}|_{\Lambda_{\mathrm{opt}}}\leq s}V(\tilde{\lambda})$. It can be seen that $\phi_1$ and $\phi_2$ are positive definite and non-decreasing. Since $|\tilde{\lambda}|_{\Lambda_{\mathrm{opt}}}$ is a continuous function of $\tilde{\lambda}$ for any given $\Lambda_{\mathrm{opt}}$, the functions $\phi_1$ and $\phi_2$ are continuous. Because $V(\tilde{\lambda})$ is radially unbounded and $\Lambda_{\mathrm{opt}}$ is compact, $\phi_1(s)$ and $\phi_2(s)$ tend to infinity as $s\rightarrow\infty$. Hence we can choose $\psi_1(s)\leq k_1\phi_1(s)$ with $k_1\in(0,1)$ and $\psi_2(s)\geq k_2\phi_2(s)$ with $k_2\in(1,\infty)$ to belong to class $\mathcal{K}_{\infty}$. Therefore, $\psi_1(|\tilde{\lambda}|_{\Lambda_{\mathrm{opt}}})\leq\phi_1(|\tilde{\lambda}|_{\Lambda_{\mathrm{opt}}})\leq V(\tilde{\lambda})\leq\phi_2(|\tilde{\lambda}|_{\Lambda_{\mathrm{opt}}})\leq\psi_2(|\tilde{\lambda}|_{\Lambda_{\mathrm{opt}}})$ holds. By Rademacher's theorem, $V(\tilde{\lambda})$ is differentiable almost everywhere and $\langle \frac{\partial V}{\partial \tilde{\lambda}},-\Theta(t)\nabla f(\tilde{\lambda})\rangle=-\nabla^Tf(\tilde{\lambda})\Theta(t)\nabla f(\tilde{\lambda})$ holds for almost all $\tilde{\lambda}$. Thus, from \cite[item 5, p. 100]{TP}, for almost all $t$ and each $v\in-\Theta(t)\partial f(\tilde{\lambda})$, $\dot{V}(\tilde{\lambda})\leq V^{\circ}(\tilde{\lambda};v)\leq-\langle\nabla f(\tilde{\lambda}),\Theta(t)\nabla f(\tilde{\lambda})\rangle$. Since $\Theta(t)$ is diagonal and all its diagonal elements are in $[\varepsilon,1]$, $\nabla^Tf(\tilde{\lambda})\Theta(t)\nabla f(\tilde{\lambda})$ is lower bounded by the continuous function $\varepsilon\nabla^Tf(\tilde{\lambda})\nabla f(\tilde{\lambda})$ which vanishes only on $\Lambda_{\mathrm{opt}}$ and is positive everywhere else. Thus there exists a positive definite function $\psi_3(s):=\inf_{|\tilde{\lambda}|_{\Lambda_{\mathrm{opt}}}\geq s}\varepsilon\nabla^Tf(\tilde{\lambda})\nabla f(\tilde{\lambda})$, such that $\dot{V}\leq-\psi_3(|\tilde{\lambda}|_{\Lambda_{\mathrm{opt}}})$ for almost all $t$. By \cite[Lemma IV.1]{ASW}, there exist $\rho_1\in\mathcal{K}_{\infty}$ and $\rho_2\in\mathcal{L}$ such that $\psi_3(s)\geq\rho_1(s)\rho_2(s)$. Therefore, for almost all $t$, $\dot{V}\leq-\psi_3(|\tilde{\lambda}|_{\Lambda_{\mathrm{opt}}})\leq-\rho_1(|\tilde{\lambda}|_{\Lambda_{\mathrm{opt}}})\rho_2(|\tilde{\lambda}|_{\Lambda_{\mathrm{opt}}})\leq-\rho_1(\psi_2^{-1}(V))\rho_2(\psi_1^{-1}(V)):=-\rho(V)$. Since $\rho$ is positive definite, by \cite[Lemma IV.2]{ASW}, there exists a $\beta\in\mathcal{KL}$ such that $V(\tilde{\lambda}(t))\leq\beta(V(\tilde{\lambda}(t_0)),t-t_0)$. Therefore, all solutions to (\ref{e2.1}) starting at $\tilde{\lambda}(t_0)$ satisfy $|\tilde{\lambda}(t)|_{\Lambda_{\mathrm{opt}}}\leq \psi_1^{-1}(\beta(\psi_2(|\tilde{\lambda}(t_0)|_{\Lambda_{\mathrm{opt}}}),t-t_0)):=\bar{\beta}(|\tilde{\lambda}(t_0)|_{\Lambda_{\mathrm{opt}}},t-t_0)$ with $\bar{\beta}\in\mathcal{KL}$. Thus, $\Lambda_{\mathrm{opt}}$ is UGAS for (\ref{e2.1}).

%\begin{IEEEbiography}[{\includegraphics[width=1in,height=1.25in,clip,keepaspectratio]{Yankai.jpg}}]{Yankai Lin} (S'18) is currently pursuing the Ph.D. degree at the University of Melbourne. He excels at writing concise biographies. 
%\end{IEEEbiography}
%
%\begin{IEEEbiography}[{\includegraphics[width=1in,height=1.25in,clip,keepaspectratio]{a2.jpg}}]{Iman Shames} (M'11) is an Associate Professor and a squash player at The University of Melbourne. 
%\end{IEEEbiography}
%
%\begin{IEEEbiography}[{\includegraphics[width=1in,height=1.25in,clip,keepaspectratio]{a3.png}}]{Dragan Ne\v{s}i\'{c}} (F'08) is a Melbourne-based guitarist who proudly runs a private gym and an olive farm at the same time. 
%\end{IEEEbiography}

\end{document}